\documentclass[a4paper]{article}%

\usepackage{amsmath,amsthm,amssymb,fancyhdr,color,epic,graphicx}%
\usepackage[titletoc,toc,title]{appendix}

\newtheorem{theorem}{\bf Theorem}[section]

\newtheorem{lemma}[theorem]{\bf Lemma}
\newtheorem{proposition}[theorem]{\bf Proposition}

\newenvironment{keywords}{\par\noindent\textbf{Keywords: }}{\hfill\par}

\theoremstyle{definition}

\newtheorem{remark}[theorem]{\bf Remark}

\definecolor{Be}{rgb} {0,0.08,0.45} 
\definecolor{darkgreen}{rgb}{0.00, 0.49, 0.00}
\definecolor{Mahog}{rgb}{0.6,  0.,   0}

\newcommand{\R}{I\!\!R}
\newcommand{\C}{\mathbb{C}}
\newcommand{\N}{I\!\!N}

\numberwithin{equation}{section}

\begin{document}

\title{The Landau-Zener transition and the surface hopping method for the 2D Dirac equation for graphene\footnote{Research supported by 
NSFC grant 91330203, NSF grant 1114546 and  the NSF grant RNMS-1107291 
``KI-Net''}}
\author{Ali Faraj\footnote{Institute of Natural Sciences, Shanghai Jiao Tong University, Shanghai 200240 P.R. China. Current address: Grenoble INP, ESISAR, 26902 Valence Cedex 9, France. Email: ali.faraj@esisar.grenoble-inp.fr} and Shi Jin\footnote{Department of Mathematics, 
Institute of Natural Sciences, MOE-LSEC and SHL-MAC, Shanghai Jiao Tong University, Shanghai 200240 P.R. China; and Department of Mathematics, University of Wisconsin, Madison, WI 53706, USA. Email: jin@math.wisc.edu}}
\date{}

\maketitle

\begin{abstract} A Lagrangian surface hopping algorithm is implemented
to study the two dimensional massless Dirac equation for Graphene with an electrostatic potential, in the semiclassical regime. In this problem, the crossing of the energy levels of the system at Dirac points requires a particular treatment in the algorithm in order to describe the quantum transition--characterized by the Landau-Zener probability-- between
different energy levels. 
We first derive the Landau-Zener probability for the underlying problem,
then incorporate it into the surface hopping algorithm. We also show that different
 asymptotic models for this problem derived in \cite{MoSu2}
 may give  different transition
probabilities. We conduct numerical experiments to compare the solutions to
the Dirac equation, the surface hopping algorithm, and the asymptotic 
models of \cite{MoSu2}.

\bigskip

\bigskip

\end{abstract}

\section{Introduction}

We are interested in the description of the transport of electrons in a single graphene layer. This material is a two-dimensional flat monolayer of carbon atoms which displays unusual and interesting electronic properties arising from the bi-conically shaped Fermi surfaces near the Brillouin zone corners (Dirac points). The electrons propagate as massless Dirac Fermions moving with the Fermi velocity $v_F$, which is $300$ times smaller than the speed of light $v_F\approx \frac{c}{300}\approx10^6\,m.s^{-1}$, and their behavior reproduces the physics of quantum electrodynamics but at much smaller energy scale. Although this model has been studied for a long time, see \cite{CNGPNG} for a bibliography, it has remained theoretical until the work of \cite{Nall} where the graphene was produced for the first time. After these results, the interest of researchers on this material has shown a remarkable increase including applications in carbon-based electronic devices \cite{LEBK} and numerical simulations, see e.g. \cite{FiIa} and references therein.

In this paper, we will consider a  model of a two-dimensional Dirac equation (\cite{Bee, CNGPNG, No}) of a graphene sheet in the presence of an external potential. This model consists of a small parameter $h$ directly related to the Planck constant. 
We are interested is the design of an efficient numerical method--the surface hopping method-- for the graphene Dirac equation in the semiclassical regime
where $h<<1$. In this regime,
 the solution of the Dirac equation is  highly oscillatory thus a huge computational cost is required to give accurate wave functions or physical observables
for either finite difference methods (\cite{BrHeMa}, \cite{HPA}) or time-splitting spectral method \cite{HJMSZ}, since one needs to resolve the high  frequency
both spatially and temporally.  

The development of efficient numerical methods for the related Schr\"odinger equation in the semiclassical regime has motivated many works in the last decade, see the review paper \cite{JMS} and references therein. In the semiclassical regime, one often  uses  asymptotic 
analysis, such as the WKB analysis and the Wigner transform to help to reduce the computational
costs and to develop efficient computational methods based on the asymptotic
models. In the framework of Wigner transform, the idea is to construct a measure on the phase space, called the Wigner measure,  when $h\rightarrow0$, 
to obtain the
physical observables (such as density, flux, and energy) with a computational
cost far less than a direct quantum simulation. When the gap between different
 energy levels is of order one (the so-called {\it adiabatic} case), the Wigner measure technique provides a simple description of the motion: it can be well approximated by a fully diagonalized system, one classical Liouville equation for each energy level \cite{GMMP}. However, in the graphene Dirac equation, the energy levels {\it cross} at the Dirac points, where {\it non-adiabatic} transfers are observed and the particles can tunnel from one band to the other. The standard Wigner approach then needs to be revised to describe
the non-adiabatic phenomena. 

 One of the widely used approaches to simulate the non-adiabatic dynamics is the surface hopping method initially proposed by Tully and Preston \cite{TuPr}
as an efficient computational method to go beyond the classical Born-Oppenheimer
approximation. This method is widely used in chemistry and molecular 
dynamics, see for examples
\cite{Dru,ST, Tul,  TuPr}.
 The basic idea is to combine the classical transports of the system on individual potential energy surfaces with instantaneous transitions from one energy surface to the other. The transition rates for this band-to-band hopping are given by the well known Landau-Zener formula \cite{Ze}. From the mathematical point of view, the first rigorous analysis of non-adiabatic transfer using Landau-Zener formula dates back to Hagedorn \cite{Ha}. More recently, the Wigner measure techniques for separated energy levels have been extended in \cite{FeGe1} and \cite{FeGe2} to systems presenting band crossing. The proof is based on microlocal analysis to justify the Landau-Zener formula. Departing from the results in \cite{FeGe1} and \cite{FeGe2}, a rigorous surface hopping algorithm in the Wigner picture was proposed in \cite{LaTe} for time-dependent two-level Schr\"odinger systems with conically intersecting eigenvalues and implemented numerically in the Lagrangian formulation in \cite{LaSwTe}. The corresponding Eulerian numerical scheme was proposed in \cite{JiQiZa} by formulating the hopping mechanism as an interface condition which is then built into the numerical flux for solving the underlying Liouville equation for each energy level.

 In the present article we give a Lagrangian surface hopping algorithm 
for the graphene Dirac equation similar to the algorithm in \cite{LaSwTe}. 
First, it is a classical result that the Wigner transform leads to two decoupled classical Liouville equations for each energy level, under the adiabatic assumption \cite{GMMP}.  At the Dirac points where non-adiabatic transition occurs, we first derive the  Landau-Zener formula, and then incorporate it
into the surface hopping algorithm.  We also show that a reduced asymptotic
model developed in \cite{MoSu2} could give incorrect transition probability.
We then compare through several numerical examples the solutions of the Dirac equation
(solved by the time-splitting spectral method),
the surface hopping algorithm and the asymptotic models of \cite{MoSu2}.
Our numerical results  show that, when there is no wave interference,
 the surface hopping algorithm indeed
gives the correct non-adiabatic transition at Dirac points, with a
much greater computational efficiency compared with the simulations based on
 solving
directly the Dirac equation. 


 The article is organized as follows: in section \ref{sec_pb} we give  the graphene Dirac equation and its semiclassical limit via the Wigner transform
in the adiabatic case. We give some examples of potentials to show that
non-adiabatic transition is indeed possible.  In section \ref{sec_non_ad_tr}, 
we derive the Landau-Zener transition probability for the graphene Dirac equation. The surface hopping algorithm is given in section \ref{sec_alg}.  In section \ref{sec_thmo}, we study the asymptotic models introduced  in \cite{MoSu2} and
show that a reduced model could give the incorrect transition probability.
Numerical results are given in section \ref{sec_numres} for comparisons of
different models. For reader's convenience we give the time-splitting
spectral method of the Dirac equation in Appendix \ref{sec_Dsol}. 

\section{Quantum transport in graphene and the Wigner measure}\label{sec_pb}

We consider the description of the transport of electrons in a single graphene layer in the presence of an external potential. Following \cite{Bee, CNGPNG, No}, it is modelled by the two dimensional Dirac equation:
\begin{equation}\label{eq_dirsyst}
\displaystyle \left\{\begin{array}{ll} i\hbar\partial_t\psi = \left[-i\hbar v_F\sigma_1\partial_{x_1}-i\hbar v_F\sigma_2\partial_{x_2}+qV\right]\psi, & t\in\R, \, x\in\R^2 \\
\psi(0,x)=\psi_I(x),&x\in\R^2
\end{array}\right.
\end{equation}
where $\hbar$ is the reduced Planck constant, $v_F$ is the Fermi velocity, $q$ is the elementary charge of the electron and $V(x)\in\R$ is the electric potential. The Pauli matrices $\sigma_1,\,\sigma_2$ are given by:
\[
\sigma_1 = \left(\begin{array}{cc} 0 & 1 \\
1 & 0 \end{array}\right), 
\quad 
\sigma_2 =\left(\begin{array}{cc} 0 & -i \\
i & 0 \end{array}\right)\,.
\]
The initial wave function $\psi_I(x)\in\C^2$ is normalized such that:
\begin{equation}\label{eq_norm1}
\int_{\R^2}|\psi_I(x)|^2dx = 1
\end{equation}
and, using the mass conservation, the wave function $\psi(t,x)\in\C^2$ 
satisfies
\begin{equation}\label{eq_norm2}
\int_{\R^2}|\psi(t,x)|^2dx = 1
\end{equation}
where t, $x=(x_1,x_2)$ are the time and space variables respectively. We  consider the system \eqref{eq_dirsyst} in the semiclassical regime. For this purpose, we rewrite the equations such that there remains only one dimensionless parameter $h$. Proceeding as in \cite{SpMa}, we change the variables to dimensionless variables as follows
\begin{equation}\label{eq_changdev0}
t \rightarrow t/T_0, \quad x \rightarrow x/L 
\end{equation}
and define
\begin{equation}\label{eq_changdev}
u_I(x) = L \psi_I(Lx), \quad u(t,x) = L \psi(T_0t,Lx)
\end{equation}
where $L$ is a reference length and $T_0=L/v_F$. We remark that $u_I(x)$ and $u(t,x)$ are chosen such that the change of variable preserves the normalization \eqref{eq_norm1} and \eqref{eq_norm2}. Plugging \eqref{eq_changdev} into \eqref{eq_dirsyst}, and dividing by the reference energy $mv_F^2$ where $m$ is the effective mass of the electron, one gets the dimensionless graphene Dirac equation:
\begin{equation}\label{eq_dirad}
\displaystyle \left\{\begin{array}{ll} ih\partial_tu = \left[-ih\sigma_1\partial_{x_1}-ih\sigma_2\partial_{x_2}+U\right]u, & t\in\R, \, x\in\R^2 \\
u(0,x)=u_I(x),&x\in\R^2
\end{array}\right.
\end{equation}
where
\begin{equation}\label{eq_defh}
h=\frac{\hbar}{mv_FL}
\end{equation}
is a small dimensionless parameter and $U(x)$ is the dimensionless potential defined by
\[
U(x) = \frac{q}{mv_F^2}V(Lx)\,.
\] 
We will consider the solution of \eqref{eq_dirad} in the limit $h\rightarrow 0$.

Non-adiabatic transfer happens at the Dirac points which are the crossing points of the eigenvalues of the symbol related to \eqref{eq_dirad}. To be more precise, we consider the complex $2\times 2$-matrix-valued symbol:
\begin{equation}\label{eq_sP0}
P_0(x,\xi) = B(\xi) + U(x)I, \quad (x,\xi)\in \R^2_x\times\R^2_\xi\,,
\end{equation}
where $I$ is the $2\times2$ identity matrix, $B(\xi)$ is given by
\begin{equation*}
B(\xi) = \xi_1\sigma_1+\xi_2\sigma_2=\left(\begin{array}{cc} 0 & \xi_1-i\xi_2 \\
\xi_1+i\xi_2& 0
\end{array}\right), \quad \xi\in\R^2,
\end{equation*}
and $U\in C^{\infty}(\R^2,\R)$ is such that $\forall\beta\in\N_0^2$ there is a constant $C_{\beta}>0$ verifying:
\[
\left|\partial_x^{\beta}U(x)\right|\leq C_{\beta}, \quad \forall x\in \R^2\,.
\]
An easy computation shows that the matrix $B(\xi)$ has eigenvalues $\pm|\xi|$ with corresponding orthonormal set of eigenvectors given by
\begin{equation}\label{eq_vectprop}
\chi_{\pm}(\xi) = \frac{1}{\sqrt{2}}\left(1,\pm\frac{\xi_1+i\xi_2}{|\xi|}\right)^T\,.
\end{equation}
Moreover, it holds $P_0\in S^1(\R^2\times\R^2)^{2\times2}$ and the Dirac equation \eqref{eq_dirad} can be rewritten as:
\begin{equation}\label{eq_dir}
\displaystyle \left\{\begin{array}{ll} ih\partial_tu^h = P_0(x,hD)u^h, & t\in\R \\
u^h(0)=u^h_I
\end{array}\right.
\end{equation}
where $P_0(x,hD)$ is the Weyl operator defined for $u\in\C^{\infty}_0(\R^2)^2$ by the integral:
\begin{equation}\label{eq_wop}
P_0(x,hD)u(x)=\frac{1}{(2\pi)^2}\int_{\R^2_\xi}\int_{\R^2_y}P_0\left(\frac{x+y}{2},h\xi\right)u(y)e^{i(x-y).\xi}d\xi dy\,,
\end{equation}
and $D=-i\partial_x$. The operator $P_0(x,hD)$ is essentially self-adjoint on the Hilbert space $\mathcal{H}=L^2(\R^2)^2$ and the domain of its self-adjoint extension is $H^1(\R^2)^2$. Therefore $-\frac{i}{h}P_0(x,hD)$ generates a strongly continuous group of unitary operators solution to \eqref{eq_dir}.

The eigenvalues $\lambda_{\pm}(x,\xi)= U(x)\pm \vert\xi\vert$ of the symbol $P_0(x,\xi)$ satisfy $\lambda_{+}(x,\xi)=\lambda_{-}(x,\xi)$ at the crossing set $\{\xi=0\}\subset\R^2_x\times\R^2_\xi$. The semiclassical limit away from the crossing set was performed in \cite{GMMP} for systems of the form \eqref{eq_dir} with more general symbols and initial data $u^h_I$ subject to additional conditions. In \cite{GMMP}, the authors show that the semiclassical limit for the different bands can be treated separately where, for each band, the related eigenvalue plays the role of a scalar classical Hamiltonian.\\

\subsection{Adiabatic semiclassical limit}\label{adiabatic}

We recall now some basic notions of the Wigner analysis involved in the semiclassical limit. For $u,\, v\in L^2(\R^2)$, the Wigner transform is defined by:
\[
w^h(u,v)(x,\xi)= \frac{1}{(2\pi)^2}\int_{\R^2}u\left(x-h\frac{y}{2}\right)\overline{v}\left(x+h\frac{y}{2}\right)e^{i\xi.y}dy
\]
and for $u\in \mathcal{H}$, the $2\times 2$ Wigner matrix is defined by:
\[
W^h[u] = \left( w^h(u_i,u_j) \right)_{1\leq i,j \leq 2}\,.
\]
We denote by $w^h[u]=\textrm{tr}\,W^h[u]$ the scalar Wigner transform of $u$. For any bounded sequence $f^h$ in $\mathcal{H}$, there is a subsequence of $W^h[f^h]$ which converges in $\mathcal{S}'$. Such a limit $W^0$ is called a Wigner measure associated to $f^h$. If $f^h$ admits only one Wigner measure, we shall denote it by $W^0[f^h]$ and  set $w^0[f^h] = \textrm{tr}\,W^0[f^h]$.

Another important object is the classical flow $\phi^{\pm}_t(x,\xi)=(x^{\pm}(t),\xi^{\pm}(t))$ corresponding to the eigenvalues $\lambda_{\pm}$, i.e. the solution to:
\begin{equation}\label{eq_caract}
\displaystyle \left\{\begin{array}{ll}\frac{d}{dt}x^{\pm}(t)=\pm \frac{\xi^{\pm}(t)}{\vert\xi^{\pm}(t)\vert}, & x^{\pm}(0)=x\\
\frac{d}{dt}\xi^{\pm}(t)=-\partial_xU\left(x^{\pm}(t)\right), & \xi^{\pm}(0)=\xi
\end{array}\right.\,,
\end{equation} 
where $(x,\xi)\in\R^2_x\times\R^2_\xi$. Indeed, the decoupled semiclassical limit is valid on a set of the phase space which is stable under the flow $\phi^{\pm}_t$ and where no band crossing occurs. More precisely, if there exists an open subset $\Omega\subset\R^2_x\times\R^2_\xi$ such that:
\[
\Omega\cap\{\xi=0\}=\emptyset \quad \textrm{and} \quad \phi^{\pm}_t(\Omega)\subset\Omega, \; \forall t \in\R
\]
and if the initial condition $u^h_I$ in \eqref{eq_dir} has a Wigner measure $W_I^0$ such that $w_I^0 = \textrm{tr}\,W_I^0$ satisfies
\[
w_I^0|_{\Omega^c}=0
\] 
then using the results in \cite{GMMP}, it holds that $w^h[u^h]$, the scalar Wigner transform of the solution $u^h$ to \eqref{eq_dir}, converges to
\begin{equation}\label{eq_decmes}
w^0(t,x,\xi) = w^0_+(t,x,\xi) + w^0_-(t,x,\xi)
\end{equation}
where $w^0_{\pm}(t,.,.)$ is the scalar positive measure on $\R^2_x\times\R^2_\xi$ solving the classical Liouville equations:
\begin{equation}\label{eq_liouv}
\displaystyle \left\{\begin{array}{ll} \partial_t w^0_{\pm} \pm \frac{\xi}{\vert\xi\vert}\partial_xw^0_{\pm}-\partial_xU.\partial_{\xi}w^0_{\pm}=0, \quad \R_t\times\Omega\\
w^0_{\pm}(0,.,.) = \textrm{tr}\,\left(\Pi_{\pm}W_I^0\right), \quad \Omega & w^0_{\pm}(t,\Omega^c) = 0, \quad t\in\R\end{array}\right.\,.
\end{equation}
For all $(x,\xi)\in\Omega$, $\Pi_{\pm}(\xi)$ in \eqref{eq_liouv} denotes the projection on the eigenspace related to the eigenvalue $\lambda_{\pm}(x,\xi)$. Using  formula \eqref{eq_vectprop}, it can be expressed explicitly as follows:
\begin{equation}\label{eq_defproj}
\Pi_{\pm}(\xi) = \frac{1}{2}\left( I \pm \frac{1}{|\xi|}B(\xi)  \right)\,.
\end{equation}
Moreover, the density $n^h(t,x)=|u^h(t,x)|^2$ converges to
\[
n^0(t,x) = \int_{\R^2}w^0(t,x,\xi)d\xi\,.
\]
It follows that for initial data $u^h_I$ such that the bands are separated initially, i.e. the support of $w_I^0$ does not intersect with the crossing set $\{\xi=0\}$, the measure $w^0$ is described by $w_I^0$ if the support of $w_I^0$ is stable under the characteristics solution to \eqref{eq_caract} (in that case, the characteristics starting from the support of $w_I^0$ do not reach $\{\xi=0\}$).\\

\subsection{Some case studies}\label{sec_part}

\subsubsection{Case $U=0$}\label{seq_U0}

In the case of the trivial potential, the solutions to \eqref{eq_caract} are given by:
\begin{equation}\label{eq_carU0}
\left(x^{\pm}(t),\xi^{\pm}(t)\right)=\left(x\pm \frac{\xi}{|\xi|}t,\xi\right)
\end{equation}
for $\xi\neq 0$ and we can take $\Omega=\R^2_x\times\R^2_\xi\setminus\{\xi=0\}$. Then the system \eqref{eq_liouv} writes:
\begin{equation}\label{eq_liouvU0}
\displaystyle \left\{\begin{array}{ll} \partial_t w^0_{\pm} \pm \frac{\xi}{\vert\xi\vert}\partial_xw^0_{\pm}=0, \quad \R_t\times\{\xi\neq 0\}\\
w^0_{\pm}(0,.,.) =w_{I,\pm}^0, \quad \{\xi\neq 0\}& w^0_{\pm}(t,x,0) = 0, \quad t\in\R,\,x\in\R^2\end{array}\right.
\end{equation}
where
\[
w_{I,\pm}^0=\textrm{tr}\,\left(\Pi_{\pm}W_I^0\right), \quad \xi\neq 0\,.
\]
Using \eqref{eq_carU0} in \eqref{eq_liouvU0}, we get that:
\[
w^0_{\pm}(t,x,\xi) = w_{I,\pm}^0\left(x\mp\frac{\xi}{|\xi|}t,\xi\right), \quad \xi\neq 0\,.
\]
Therefore, the density $n^h(t,x)=|u^h(t,x)|^2$ converges to
\[
n^0(t,x) = \int_{\R^2}w_{I,-}^0\left(x+\frac{\xi}{|\xi|}t,\xi\right)d\xi + \int_{\R^2}w_{I,+}^0\left(x-\frac{\xi}{|\xi|}t,\xi\right)d\xi\,.
\]
In the present case, if the support of the initial scalar Wigner measure $w_I^0$ does not contain the point $\{\xi=0\}$, then no hopping occurs: the bands do not communicate and at any time the measure $w^0$ is described by the separate evolution of the level characteristics.  

\subsubsection{Case $U=\alpha x_1$, $\alpha\in\R\setminus\{0\}$}\label{seq_Uax1}

In that case, the solutions to \eqref{eq_caract} are such that:
\begin{equation}\label{eq_xit}
\xi^{\pm}(t) = \xi - (\alpha t, 0)
\end{equation}
and we can take $\Omega=\R^2_x\times\R^2_\xi\setminus\{\xi_2=0\}$. Then, system \eqref{eq_liouv} becomes:
\begin{equation}\label{eq_liouvUax1}
\displaystyle \left\{\begin{array}{l} \partial_t w^0_{\pm} \pm \frac{\xi}{\vert\xi\vert}\partial_xw^0_{\pm}-\alpha\partial_{\xi_1}w^0_{\pm}=0, \quad \R_t\times\{\xi_2\neq 0\}\,,\\
w^0_{\pm}(0,.,.) = \textrm{tr}\,\left(\Pi_{\pm}W_I^0\right), \quad \{\xi_2\neq 0\} \,, \quad\quad \, w^0_{\pm}(t,x,(\xi_1,0)) = 0, \quad t\in\R,\,x\in\R^2,\,\xi_1\in\R\end{array}\right.\,.
\end{equation}
Unlike in section \ref{seq_U0}, non-adiabatic transfer may occur. Indeed, if the support of $w_I^0$ does not intersect with the crossing set $\{\xi=0\}$ but contains points of the form $(x,\xi)$ where $\xi=(\xi_1,0)$ and $\xi_1\neq 0$, then the characteristic curve $\xi^{\pm}(t)$, given by \eqref{eq_xit}, starting from $(\xi_1,0)$, will reach the crossing set at the time $t=\frac{\xi_1}{\alpha}$ and a band-to-band transition will take place.\\
In these conditions, the system \eqref{eq_liouvUax1} does not describe correctly the asymptotics $h\rightarrow 0$ of the Wigner matrix $W^h[u^h]$. 

The goal of this paper is to develop efficient semiclassical methods to
compute the non-adiabatic transition between different bands.
To quantify the transition rates, we propose in section \ref{sec_non_ad_tr} 
a Landau-Zener formula which can be justified theoretically by using the results in \cite{FeGe2}.
Indeed,  in \cite{FeGe2}  the following symbol was considered:
\begin{equation}\label{eq_sP1}
P_1(x,\xi) = A(\xi) + U(x)I, \quad (x,\xi)\in \R^2_x\times\R^2_\xi\,,
\end{equation}
where $I$ is the $2\times2$ identity matrix and
\begin{equation}\label{eq_matA}
A(\xi) = \left(\begin{array}{cc} \xi_1 & \xi_2 \\
\xi_2& -\xi_1
\end{array}\right)\,.
\end{equation}   
Using the unitary equivalence
\begin{equation}\label{eq_equn}
B(\xi) = R^*A(\xi)R\,,
\end{equation}
where
\begin{equation}\label{eq_matR}
R = \frac{1}{\sqrt{2}}\left(\begin{array}{cc} 1 & 1 \\
i&-i
\end{array}\right)\,,
\end{equation}
it follows that the two symbols have the same eigenvalues and $P_0$ can be reduced to $P_1$ after conjugation with the matrix $R$ in the treatment of band crossing.
\begin{remark}
The functions $w^0_{\pm}(t,x,\xi)$ in \eqref{eq_decmes} are the diagonal terms of the Wigner measure $W^0(t,x,\xi):=W^0[u^h(t)](x,\xi)$. Indeed, it was shown in \cite{GMMP} that $w^0_{\pm}$ is given by
\begin{equation}\label{eq_w0}
w^0_{\pm} = \textrm{tr}\,\left(\Pi_{\pm}W^0\right)1_{\Omega}
\end{equation}
and that $W^0$ is diagonal in the sense that $W^0 = \Pi_+W^0\Pi_+ + \Pi_-W^0\Pi_-$ on $\R_t\times\Omega$. Moreover, an easy computation leads to:
\begin{equation}\label{eq_diag}
\Pi_{\pm}W^0\Pi_{\pm} = \textrm{tr}\,\left(\Pi_{\pm}W^0\right)\Pi_{\pm} = w^0_{\pm}\Pi_{\pm} \quad \textrm{ on } \quad \R_t\times\Omega\,.
\end{equation} 
\end{remark}

\section{The Landau-Zener transition}\label{sec_non_ad_tr}

In this section we will study the  Landau-Zener transition for the Hamiltonian $P_0(x,hD)$.

\subsection{Classical flow around the crossing set}

As remarked in section \ref{sec_pb}, non-adiabatic transfer happens only when the characteristics reach the crossing set $\{\xi=0\}$. Proposition \ref{prop_cross} below says that such characteristics will exist as soon as the potential has points $x\in\R^2$ such that $\partial_xU(x)\neq0$. We refer to \cite{FeGe2} for the proof.
\begin{proposition}\label{prop_cross}
Consider $x\in\R^2$ such that $\partial_xU(x)\neq0$, then there exist two unique curves $s\mapsto (x^{\pm}(s),\xi^{\pm}(s))$ which are continuous in $s$ in a neighborhood of $0$ and $C^1$ for $s\neq0$ and such that:
\begin{equation}\label{eq_car0}
\displaystyle \left\{\begin{array}{ll}\frac{d}{ds}x^{\pm}(s)=\pm \frac{\xi^{\pm}(s)}{\vert\xi^{\pm}(s)\vert}, & x^{\pm}(0)=x\\
\frac{d}{ds}\xi^{\pm}(s)=-\partial_xU\left(x^{\pm}(s)\right), & \xi^{\pm}(0)=0
\end{array}\right.\,.
\end{equation}
\end{proposition}
To illustrate the fact that a spectral transfer for the symbol $P_0$ happens at the crossing set, we will come back momentarily to the potential $U=\alpha x_1$ introduced in section \ref{seq_Uax1}. For such a potential, problem \eqref{eq_car0} writes:
\[
\displaystyle \left\{\begin{array}{ll}\frac{d}{ds}x^{\pm}(s)=\pm \frac{\xi^{\pm}(s)}{\vert\xi^{\pm}(s)\vert}, & x^{\pm}(0)=x\\
\frac{d}{ds}\xi^{\pm}(s)=-\alpha\left(\begin{array}{c}1\\0\end{array}\right), & \xi^{\pm}(0)=0
\end{array}\right.\,,
\] 
and $\forall x\in\R^2$, its unique solution is given by:
\[
x^{\pm}(s) = x \mp \textrm{sgn}(\alpha)\left(\begin{array}{c}|s|\\0\end{array}\right), \quad \xi^{\pm}(s)=-\alpha\left(\begin{array}{c}s\\0\end{array}\right)\,.
\]
Plugging this solution in the projectors defined in \eqref{eq_defproj}, we get for $s<0$:
\[
\Pi_{+}(\xi^{\pm}(s)) = \frac{1}{2}\left(\begin{array}{cc} 1 & \textrm{sgn}(\alpha) \\ \textrm{sgn}(\alpha) & 1 \end{array}\right), \quad 
\Pi_{-}(\xi^{\pm}(s)) = \frac{1}{2}\left(\begin{array}{cc} 1 & -\textrm{sgn}(\alpha) \\ -\textrm{sgn}(\alpha) & 1 \end{array}\right)
\]
and for $s>0$:
\[
\Pi_{+}(\xi^{\pm}(s)) = \frac{1}{2}\left(\begin{array}{cc} 1 & -\textrm{sgn}(\alpha) \\ -\textrm{sgn}(\alpha) & 1 \end{array}\right), \quad 
\Pi_{-}(\xi^{\pm}(s)) = \frac{1}{2}\left(\begin{array}{cc} 1 & \textrm{sgn}(\alpha) \\ \textrm{sgn}(\alpha) & 1 \end{array}\right)\,.
\]
It is easy to see that the projectors $\Pi_{+}(\xi^{\pm}(s))$ and $\Pi_{-}(\xi^{\pm}(s))$ interchange one with the other when the characteristics pass through the crossing set $\{\xi=0\}$.
 
\subsection{A heuristic derivation of the Landau-Zener formula}\label{sec_heur}

In this section, we give a heuristic argument, similar to the one in \cite{LaSwTe} and \cite{LaTe}, to derive the Landau-Zener formula for the Hamiltonian $P_0(x,hD)$.

In general,  the region for non-adiabatic transfer is not restricted to the crossing set $\{\xi=0\}$, since quantum transition occurs as long as the two
energy levels are sufficiently close, in the case of avoided crossing \cite{Ha}.  As in \cite{LaSwTe} and \cite{LaTe},
we define this region to contain the points where the distance between the eigenvalues $\lambda_{\pm}(x,\xi)$ of $P_0(x,\xi)$ is minimal. In our case $|\lambda_+(x,\xi)-\lambda_-(x,\xi)|=2|\xi|$ and, when considered along the characteristics solution to \eqref{eq_caract}, the necessary condition for minimal gap is:
\[
\left(\vert\xi(s)\vert^2\right)'=0 \Leftrightarrow \xi(s).\partial_xU(x(s))=0
\]
and the hopping surface is chosen as the set:
\begin{equation}\label{eq_minsurf}
S = \{ (x,\xi)\in\R^4;\, \xi.\partial_xU(x)=0  \}\,.
\end{equation}
The heuristics follows by inserting the characteristics in the trace-free part of the symbol to obtain the system of ordinary differential equations:
\[
ih\psi'(s) = B(\xi(s))\psi(s)\,.
\]
After conjugation with the matrix $R$ defined in \eqref{eq_matR} and using equation \eqref{eq_equn}, one arrives at the following new system:
\[
ih\psi'(s) = A(\xi(s))\psi(s)\,.
\]
Assume the particles defined by the trajectory (\eqref{eq_caract}) are near a point $(x^*,\xi^*)\in S$ (due to translation invariance we assume  particles are
at $(x^*,\xi^*)$ initially), then the Taylor expansion gives:
\[
x(s) = x^* \pm \frac{\xi^*}{\vert\xi^*\vert}s + \mathcal{O}(s^2)
\]
\[
\xi(s) = \xi^* -\partial_xU(x^*)s + \mathcal{O}(s^2)\,.
\]
Ignoring the $\mathcal{O}(s^2)$ terms, the system becomes:
\[
ih\psi'(s) = \left(\begin{array}{cc} \xi^*_1 -\partial_{x_1}U(x^*)s  & \xi^*_2 -\partial_{x_2}U(x^*)s \\
\xi^*_2 -\partial_{x_2}U(x^*)s& -\xi^*_1 +\partial_{x_1}U(x^*)s
\end{array}\right)\psi(s)\,.
\]
After conjugation with the rotation matrix:
\[
\left(\begin{array}{cc} \cos\theta & \sin\theta \\
-\sin\theta& \cos\theta
\end{array}\right)
\]
where $\theta$ is the angle such that
\[
(\cos 2\theta,\sin 2\theta) = \frac{\partial_{x}U(x^*)}{\vert\partial_{x}U(x^*)\vert}
\]
we get:
\[
i\frac{h}{\vert\partial_{x}U(x^*)\vert}\psi'(s) = \left(\begin{array}{cc} -s  & \frac{\xi^*\wedge\partial_{x}U(x^*)}{\vert\partial_{x}U(x^*)\vert^2} \\
\frac{\xi^*\wedge\partial_{x}U(x^*)}{\vert\partial_{x}U(x^*)\vert^2}& s
\end{array}\right)\psi(s)
\]
where $\xi\wedge\zeta = \xi_2\zeta_1 - \xi_1\zeta_2$ for $\xi,\,\zeta\in\R^2$. After conjugation with the matrix
\[
\left(\begin{array}{cc} 0 & 1 \\
1&0
\end{array}\right)
\]
it follows:
\[
i\frac{h}{\vert\partial_{x}U(x^*)\vert}\psi'(s) = \left(\begin{array}{cc} s  & \frac{\xi^*\wedge\partial_{x}U(x^*)}{\vert\partial_{x}U(x^*)\vert^2} \\
\frac{\xi^*\wedge\partial_{x}U(x^*)}{\vert\partial_{x}U(x^*)\vert^2}& -s
\end{array}\right)\psi(s)\,.
\]
If we set $\varepsilon = \frac{h}{\vert\partial_{x}U(x^*)\vert}$ and $\eta=\frac{\xi^*\wedge\partial_{x}U(x^*)}{\vert\partial_{x}U(x^*)\vert^2}$, the system becomes:
\[
i\varepsilon\psi'(s) = \left(\begin{array}{cc} s  & \eta \\
\eta & -s
\end{array}\right)\psi(s)
\]
which is the well known Landau-Zener problem (see \cite{Ze}) for which the transition probability is:
\[
T = e^{-\frac{\pi}{\varepsilon}\eta^2}\,.
\]
This allows us to propose the following non-adiabatic transition rate at the point $(x^*,\xi^*)$:
\begin{equation}\label{eq_tprob}
T(x^*,\xi^*) = e^{-\frac{\pi}{h}\frac{\left(\xi^*\wedge\partial_{x}U(x^*)\right)^2}{\vert\partial_{x}U(x^*)\vert^3}}\,.
\end{equation}

\subsection{About the rigorous justification of the Landau-Zener formula}\label{sec_rig}
  
As already noticed in section \ref{sec_pb}, the symbol $P_0$ involved in the Dirac equation \eqref{eq_dir} satisfies the identity
\[
P_1(x,\xi)=RP_0(x,\xi)R^*
\]
where $P_1$ and $R$ are given respectively by \eqref{eq_sP1} and \eqref{eq_matR}. Therefore, if $u^h$ denotes the solution to \eqref{eq_dir}, the function $v^h=R u^h$ satisfies the equation:
\begin{equation}\label{eq_vh2}
ih\partial_tv^h = P_1(x,hD)v^h
\end{equation}
where $P_1(x,hD)$ is defined as in \eqref{eq_wop}. A Landau-Zener formula was obtained rigorously in \cite{FeGe2} for the two-scale Wigner measure $\tilde{\nu}$ of the function $v^h$ (see \cite{FeGe2}, \cite{LaTe} for the definition of two-scale Wigner measures). This result provides a rigorous proof of our Landau-Zener formula. Indeed, if $\nu_I$ denotes the two-scale Wigner measure of $u^h$, it follows that the two-scale Wigner measure of $v^h$ is $\tilde{\nu}=R\nu_IR^*$. Then, the Landau-Zener formula for $\nu_I$ can be deduced from the Landau-Zener formula for $\tilde{\nu}$.
\begin{remark}
The Landau-Zener formula obtained in \cite{FeGe2} writes
\[
T = e^{-\frac{\pi\eta^2}{\vert\partial_{x}U(x)\vert}}
\]
where
\[
\eta = \delta\xi\wedge\frac{\partial_{x}U}{\vert\partial_{x}U\vert}\,.
\]
When the direction $\delta\xi$ is equal to $\frac{\xi}{\sqrt{h}}$, this corresponds to the Landau-Zener formula \eqref{eq_tprob} that we obtained heuristically in section \ref{sec_heur}. 
\end{remark}

\section{A surface hopping algorithm}\label{sec_alg}

We give here a semiclassical Lagrangian algorithm for the evolution of the diagonal terms of the Wigner matrix $W^h(t,x,\xi)$ of the solution $u^h(t)$ to \eqref{eq_dir}. The algorithm is adopted from the method proposed in \cite{LaSwTe} for time-dependent two-level Schr\"odinger systems with conically intersecting eigenvalues.

Define the level populations:
\begin{equation}\label{eq_lpop}
P^h_{\pm}(t) = ||\Pi_{\pm}(hD)u^h(t)||_{\mathcal{H}}^2
\end{equation}
where for $u\in\mathcal{H}$, $\Pi_{\pm}(hD)u$ is the function defined via its Fourier transform $\Pi_{\pm}(h\xi)\hat{u}(\xi)$ and $\hat{u}=\mathcal{F}u$ is given by:
\begin{equation*}
\mathcal{F}u(\xi) = \frac{1}{2\pi}\int_{\R^2}u(x)e^{-ix.\xi}dx
\end{equation*} 
(it is clear from \eqref{eq_defproj} that $\Pi_{\pm}(hD)u=\Pi_{\pm}(D)u$). With similar computations as in \cite{LaSwTe}, we obtain:
\begin{equation}\label{eq_lpopw}
P^h_{\pm}(t) = \int_{\R^2_x}\int_{\R^2_{\xi}} w_{\pm}^h(t,x,\xi) dxd\xi\,,
\end{equation}
where 
\begin{equation}\label{eq_wh}
w_{\pm}^h(t,x,\xi) = \textrm{tr}\,\left(\Pi_{\pm}(\xi)W^h(t,x,\xi)\right)\,.
\end{equation}
As it appears from the following equation
\[
\Pi_{\pm}W^h\Pi_{\pm} = w^h_{\pm}\Pi_{\pm}\,,
\]
which is obtained in a similar way as \eqref{eq_diag}, $w_{\pm}^h$ are the diagonal terms of the Wigner matrix $W^h$. Up to a small remainder, the function $w^h_{\pm}$ can be written in terms of the scalar Wigner transform of the level projections of the solution $u^h(t)$. Indeed, using Lemma 2.3 in \cite{GMMP}, it holds for all $t\in\R$:
\begin{equation}\label{eq_wigpr}
w_{\pm}^h(t) = w^h\left[u^h_{\pm}(t)\right] + o(1)
\end{equation}
in $\mathcal{D}'\left(\R^2_x\times(\R^2_\xi\setminus\{\xi=0\})\right)$ when $h\rightarrow 0$, where $u^h_{\pm}(t)$ is the function with Fourier transform given by
\[
\hat{u}^h_{\pm}(t,\xi) = \Pi_{\pm}(h\xi)\hat{u}^h(t,\xi)1_{\xi\neq0} 
\]
and $\hat{u}^h(t)$ is the Fourier transform of $u^h(t)$. Comparing the relations \eqref{eq_wh} and \eqref{eq_w0}, it follows that, in the situation of section \ref{sec_pb} without band-crossing, the partial differential equation in \eqref{eq_liouv} satisfied by $w^0_{\pm}$ can be used to approximate the time evolution of $w_{\pm}^h$. In the case of band crossing, the idea is to use \eqref{eq_liouv} for the time evolution of $w_{\pm}^h$ as long as the classical trajectories solution to \eqref{eq_caract} are away from the hopping surface $S$ defined by \eqref{eq_minsurf}. When a trajectory reaches a point $(x^*,\xi^*)\in S$ a non-adiabatic transfer of weight occurs between $w_+^h$ and $w_-^h$ with transition probability $T(x^*,\xi^*)$ given by \eqref{eq_tprob}.

\subsection*{The algorithm}

\begin{itemize}
\item[1.] Initial sampling: in this step, an appropriate sampling of the function $w_{I,\pm}^h$ defined in \eqref{eq_indiag} is chosen. Specifically, 
a set of sampling points
\[
\{(x_k,\xi_k,j_k)\in\R^2_x\times\R^2_\xi\times\{-,+\}\,; \quad k=1,...,N\}
\]
are chosen with associated weights $w_k\in\R$ given by:
\[
w_k = w_{I,j_k}^h(x_k,\xi_k)\,.
\]
\item[2.] Hopping transport: away from the set $S$, each particle $(x_k,\xi_k,j_k)$ is transported by the associated classical flow $\phi^{j_k}_t$ solution to \eqref{eq_caract}. In other words, for $t\geq 0$ small enough:
\begin{equation}\label{eq_decev}
(x_k(t),\xi_k(t)) = \phi^{j_k}_t(x_k(0),\xi_k(0)), \quad w_k(t) = w_k(0)\,. 
\end{equation}
If there exists $t^*>0$ such that $(x_k(t^*),\xi_k(t^*))=:(x_k^*,\xi_k^*)\in S$, the weight is reduced using the transition rate
\[
T^*=T\left(x_k^*,\xi_k^*\right)\,,
\] 
where $T(x,\xi)$ is given by \eqref{eq_tprob}. Moreover, in order to describe completely the non-adiabatic transfer, a new particle with index $l>N$ is created. Specifically,  for $t>t^*$
\[
(x_k(t),\xi_k(t)) = \phi^{j_k}_t(x_k(0),\xi_k(0)), \quad w_k(t) = \left(1-T^*\right)w_k(t^*) 
\]
and the new particle is created for $t>t^*$
\[
(x_l(t),\xi_l(t)) = \phi^{j_l}_{t-t^*}(x_k^*,\xi_k^*), \quad j_l = -j_k, \quad w_l(t) = T^*w_k(t^*)\,.
\]
\item[3.] Final reconstruction: at the final time $t_f>0$, there are $M\geq N$ points
\[
\{(x_k,\xi_k,j_k)\in\R^2_x\times\R^2_\xi\times\{-,+\}\,; \quad k=1,...,M\}
\]
with associated weights $w_k$ which are approximations to $w_{j_k}^h(t_f,x_k,\xi_k)$. Then, using equation \eqref{eq_lpopw}, the surface hopping approximations $P^h_{sh,\pm}(t_f)$ of the level populations $P^h_{\pm}(t_f)$ are given by:
\begin{equation}\label{eq_lpopsh}
P^h_{sh,\pm}(t_f) = \sum_{k=1}^{M}w_k\delta_{j_k}^{\pm}\omega_k
\end{equation}
where $\delta_j^i$ is the Kronecker symbol related to $i$ and $j$, and $\omega_k$ is an appropriate quadrature weight.
\end{itemize}

\begin{remark}\label{rk_noint}
We note that this surface hopping algorithm is subject to some restrictions. First only the dynamics of the diagonal components of the Wigner matrix away from the crossing set are well approximated. Second, there are possible interferences which are not captured if no particular treatment is performed. Indeed, if $w_+^h(t,x,\xi)$ and $w_-^h(t,x,\xi)$ arrive at the same time at some point $(x,\xi)$ close to the crossing set then a transfer of weight using only the Landau-Zener formula \eqref{eq_tprob} might give an incorrect approximation of the dynamics.
To avoid these interferences, we make the assumption that the initial data $u^h_I$ in \eqref{eq_dir} satisfies
\begin{equation}\label{eq_condID}
\Pi_-(h\xi)\hat{u}^h_I(\xi) = 0
\end{equation}
where $\hat{u}^h_I$ is the Fourier transform of $u^h_I$. Then, it follows from \eqref{eq_wigpr} that the diagonal terms of the initial Wigner matrix
\begin{equation}\label{eq_indiag}
w_{I,\pm}^h(x,\xi) = \textrm{tr}\,\left(\Pi_{\pm}(\xi)W^h\left[u^h_I\right](x,\xi)\right)
\end{equation}
satisfy in $\mathcal{D}'\left(\R^2_x\times(\R^2_\xi\setminus\{\xi=0\})\right)$
\begin{equation}\label{eq_dtwt}
w_{I,+}^h = w^h\left[u^h_I\right] + o(1)\,, \quad w_{I,-}^h = o(1)
\end{equation}
when $h\rightarrow 0$. Therefore, for the potential $U(x)=\alpha x_1$ given in section \ref{seq_Uax1}, the condition \eqref{eq_condID} insures that no interferences occur. Indeed, if $w_+^h(t,x,\xi)$ reaches the crossing set $\{\xi=0\}$ at some time $t^*>0$, the particle corresponding to $w_-^h(t,x,\xi)$ immediately moves away as it appears from the equation
\begin{equation}\label{eq_xitstar}
\xi^{\pm}(t) = - \alpha(t-t^*, 0)
\end{equation}
for the momentum part of the characteristics solution to \eqref{eq_caract}. Similarly, in the case where the potential has no stationary points, i.e. $\partial_xU(x)\neq0$, $\forall x\in\R^2$, equation \eqref{eq_xitstar} is replaced by 
\[
\xi^{\pm}(t) = -\int_{t^*}^t\partial_x U\left(x^{\pm}(s)\right)ds
\]
and the condition \eqref{eq_condID} insures that no interferences occur as long as $|t-t^*|$ is small enough.
\end{remark}
\begin{remark}
 To deal with this interference, a possible solution  might be to use a hybrid method for the Dirac equation \eqref{eq_dir}. Such a method was introduced in \cite{JiQi} for the Schr\"odinger equation and mixes a Gaussian beam method or a Liouville equation in the region where no-interferences occur
 and a complete quantum solver in the region where the phase information of the wave function is required. In our case, the second region corresponds to the hopping region. Since the region which involves the quantum solver can be chosen very small, a hybrid method allows an adapted treatment of the interferences without increasing dramatically the numerical cost even if the semiclassical parameter tends to $0$. In addition, the resulting algorithm is supposed to work regardless of the conditions in Remark \ref{rk_noint} on the potential and the initial data. Another possible solution is to keep the off-diagonal entries,
which contain information about the non-adiabatic transition in the Wigner transform, and then derive the semiclassical models for the entire Wigner matrix, see \cite{MoSu2, CJLM}.  In the next section, we study such a model.
 
\end{remark}

\section{Transition rate of an asymptotic model derived in \cite{MoSu2}}\label{sec_thmo}

As it was mentioned, the surface hopping method breaks down where there
are interferences. A possible remedy for this is to derive improved models
using the Wigner transform for the {\it entire} Wigner matrix, in which
the off-diagonal entries contain non-adiabatic transition information
\cite{MoSu2, CJLM}. In this section, we study such a model 
obtained in \cite{MoSu2}, and compare  the Landau-Zener transition rate of this model with the one we
derived in section \ref{sec_non_ad_tr}. To serve our purpose we will only consider the model for the potential studied in \ref{seq_Uax1}, $U=\alpha x_1$, $\alpha\in\R\setminus\{0\}$, for which the Wigner matrix does not depend on the variable $x_2$. 
It was shown in \cite{MoSu2} how this asymptotic model can describe non-adiabatic transfer as a quantum correction of the decoupled system \eqref{eq_liouvUax1}.\\ 

When the Wigner matrix does not depend on the variable $x_2$, the system \eqref{eq_liouvUax1} reads
\begin{equation}\label{eq_liouvUax1s}
\displaystyle \left\{\begin{array}{l} \partial_t w^0_{\pm} \pm \frac{\xi_1}{\vert\xi\vert}\partial_{x_1}w^0_{\pm}-\alpha\partial_{\xi_1}w^0_{\pm}=0, \quad \R_t\times\{\xi_2\neq 0\}\\
w^0_{\pm}(0,.,.) = \textrm{tr}\,\left(\Pi_{\pm}W_I^0\right), \quad \{\xi_2\neq 0\} \,; \quad\quad \, w^0_{\pm}(t,x,(\xi_1,0)) = 0, \quad t\in\R,\,x\in\R^2,\,\xi_1\in\R\end{array}\right.\,.
\end{equation}
By taking into account the scaling \eqref{eq_changdev0}\eqref{eq_changdev}, the asymptotic model in \cite{MoSu2} is the following correction of \eqref{eq_liouvUax1s}:
\begin{equation}\label{eq_liouvMo}
\displaystyle \left\{\begin{array}{l} \partial_t w_{\pm} \pm \frac{\xi_1}{\vert\xi\vert}\partial_{x_1}w_{\pm}-\alpha\partial_{\xi_1}w_{\pm}=\pm\alpha\frac{\xi_2}{|\xi|^2}\textrm{Im}\,\left(\frac{\xi_1+i\xi_2}{|\xi|}w_{i}\right), \quad \R_t\times\{\xi_2\neq 0\}\\
\partial_t w_{i} - i \Lambda(\xi) w_{i}-\alpha\partial_{\xi_1}w_{i}= -i\frac{\alpha}{2}\frac{\xi_2(\xi_1-i\xi_2)}{|\xi|^3}(w_{+}-w_{-})\end{array}\right.
\end{equation}
where 
\[
\Lambda(\xi) = -\frac{2|\xi|}{h}-\frac{\alpha\xi_2}{|\xi|^2}\,.
\]
In the limit $h\rightarrow0$, the function $w_{i}(t,x,\xi)$ approximates the off-diagonal terms of the Wigner matrix of the solution $u^h$ to \eqref{eq_dir}. As it appears from the last equation in \eqref{eq_liouvMo}, the function $w_i$ depends on $w_+$ and $w_-$. As a consequence, $w_i$ provides a coupling term at the r.h.s. of the first two equations in \eqref{eq_liouvMo}.

Using the method of characteristics, $w_i$ can be expressed as an explicit function of the difference $w_{+}-w_{-}$. Inserting this solution in the first two equations in \eqref{eq_liouvMo}, the following approximate equations for $w_{\pm}$, in the limit $h\rightarrow 0$, are obtained in \cite{MoSu2}:
\begin{equation}\label{eq_siMo}
\partial_t w_{\pm} \pm \frac{\xi_1}{\vert\xi\vert}\partial_{x_1}w_{\pm}-\alpha\partial_{\xi_1}w_{\pm}=\mp\tau(\xi)(w_{+}-w_{-}), \quad \R_t\times\R_{x_1}\times\R_{\xi_1}
\end{equation}
where $\xi_2\neq0$ can be considered as a small parameter. In the domain
\[
\left\{\xi\in\R^2\textrm{ such that }0<|\xi_2|\leq\sqrt{\alpha h}\textrm{ and }|\xi_1|\leq\frac{\alpha h}{|\xi_2|}\right\}\,,
\]
the coefficient $\tau(\xi)$ is given by:
\begin{equation}\label{eq_tauns}
\tau(\xi) = \frac{\alpha}{2}\frac{\xi_2}{|\xi|^2}\left(\frac{\pi}{2}\textrm{sgn}(\alpha\xi_2)-\arctan\frac{\xi_1}{\xi_2}\right)
\end{equation} 
where $\textrm{sgn}(x)$ denotes the sign of $x$. As it will be the case in section \ref{sec_ressh} for our surface hopping algorithm, the set $\{\xi_1=0\}$ is the important region for the hopping. Indeed, the following lemma (proof left to the reader) shows that the set $\{\xi_1=0\}$ plays the role of an interface where, in the limit $\xi_2\rightarrow0$, the solution to \eqref{eq_siMo} will have a discontinuity.
\begin{lemma}\label{le_cvtau}
The function $\xi_1 \mapsto \tau(\xi_1,\xi_2)$ tends in $\mathcal{D}'(\R)$ to $|\alpha|\beta\delta_{\xi_1=0}$ when $\xi_2 \rightarrow 0$, where
\begin{equation}\label{eq_c}
\beta = \frac{\pi^2}{4}\,.
\end{equation}
\end{lemma}  
In the case $\alpha>0$, we show below that non-adiabatic transfer is possible using the model \eqref{eq_siMo} and give the corresponding transmission matrix at the interface. Define $\omega_{\pm}$ as the set:
\[
\omega_{\pm}=(\R_t\times\R_{x_1}\times\R_{\xi_1})\cap\{\pm\xi_1>0\}
\]
and consider an initial condition in the upper band and localized at the right of the interface. In other words, the initial conditions for \eqref{eq_siMo} are
\[
w_{\pm}(0,.,.)=w_{I,\pm}
\]
where $w_{I,-}=0$ and $w_{I,+}$ is a function in $C^{\infty}_0(\R_{x_1}\times\R_{\xi_1})$ which is independent of $\xi_2$ and has support in $\omega_+$. To perform the limit $\xi_2\rightarrow 0$, the following assumptions are required on the solution to \eqref{eq_siMo}:
\begin{itemize}
\item[$A1$.] For all $\xi_2\neq0$, $w_{\pm}\in C^2(\R_t\times\R_{x_1}\times\R_{\xi_1})$.
\item[$A2$.] There is a constant $C$ which does not depend on $\xi_2$ such that:
\[
|\partial_{x_1}w_{\pm}(t,x_1,\xi_1)|\leq C, \quad \forall(t,x_1,\xi_1)\in\omega_{+}\cup\omega_{-}\,.
\] 
\item[$A3$.] There is a function $w^{0}_{\pm}$ which is continuous on $\omega_+$ and $\omega_-$ such that $w_{\pm}$ tends to $w^{0}_{\pm}$ when $\xi_2\rightarrow0$ uniformly on the compact subsets of $\omega_+$ and $\omega_-$.
\item[$A4$.] For all $t^*>0$ and $x_1^*\in\R$, the limit 
\[
\lim_{\substack{(t,x_1,\xi_1)\rightarrow(t^*,x_1^*,0)\\(t,x_1,\xi_1)\in\omega_+}}w^0_{\pm}(t,x_1,\xi_1) \quad (\textrm{resp. } \lim_{\substack{(t,x_1,\xi_1)\rightarrow(t^*,x_1^*,0)\\(t,x_1,\xi_1)\in\omega_-}}w^0_{\pm}(t,x_1,\xi_1))
\]
 exists and is denoted $w^r_{\pm}(t^*,x_1^*)$ (resp. $w^l_{\pm}(t^*,x_1^*)$).
\end{itemize}
The characteristics related to \eqref{eq_siMo} are the classical trajectories $\varphi_t^{\pm}(x_I,\xi_I)=(x_1^{\pm}(t),\xi_1^{\pm}(t))$ solution to: 
\begin{equation}\label{eq_caract2D}
\displaystyle \left\{\begin{array}{ll}\frac{d}{dt}x_1^{\pm}(t)=\pm \frac{\xi_1^{\pm}(t)}{\sqrt{\xi_1^{\pm}(t)^2+\xi_2^2}},&x_1^{\pm}(0)=x_I\\
\frac{d}{dt}\xi_1^{\pm}(t)=-\alpha,&\xi_1^{\pm}(0)=\xi_I
\end{array}\right.\,.
\end{equation}
Due to the localization of the support of the initial data $w_{I,+}$, only positive initial momenta $\xi_I>0$ have to be considered in \eqref{eq_caract2D}. Therefore, the classical trajectories will reach the interface $\{\xi_1=0\}$ at the time $t^*=\frac{\xi_I}{\alpha}$. Moreover, due to the condition $w_{I,-}=0$, the behavior at the interface is described by the characteristics corresponding to the upper band only. To be more precise, by differentiating the map $s\mapsto w_{\pm}(s,x_1^+(s),\xi_1^+(s))$, one gets by using \eqref{eq_siMo}:
\begin{equation}\label{eq_edomo}
\frac{d}{ds}\left(\begin{array}{c}w_+(s,x_1^+(s),\xi_1^+(s))\\w_-(s,x_1^+(s),\xi_1^+(s))\end{array}\right)=-\tau(\xi^+(s))M\left(\begin{array}{c}w_+(s,x_1^+(s),\xi_1^+(s))\\w_-(s,x_1^+(s),\xi_1^+(s))\end{array}\right)+F(s)
\end{equation}
where 
\[
\xi^{\pm}(s)=(\xi_1^{\pm}(s),\xi_2), \quad M=\left(\begin{array}{cc} 1 & -1\\-1 & 1\end{array}\right) \quad \textrm{ and } \quad F(s)=\frac{2\xi_1^+(s)}{|\xi^+(s)|}\left(\begin{array}{c}0\\\partial_{x_1}w_-(s,x_1^+(s),\xi_1^+(s))\end{array}\right)\,.
\] 
Applying the Duhamel formula to \eqref{eq_edomo}, it follows:
\begin{equation}\label{eq_dumo}
\left(\begin{array}{c}w_+(t,x_1^+(t),\xi_1^+(t))\\w_-(t,x_1^+(t),\xi_1^+(t))\end{array}\right)=e^{-\int_0^t\tau(\xi^+(s))dsM}\left(\begin{array}{c}w_{I,+}(x_I,\xi_I)\\w_{I,-}(x_I,\xi_I)\end{array}\right)+\int_0^te^{-\int_s^t\tau(\xi^+(\mu))d\mu M}F(s)ds\,.
\end{equation}
Using \eqref{eq_tauns}, a direct computation shows that there is a constant $C$ which does not depend on $\xi_2$ such that:
\begin{equation}\label{eq_taub}
\forall s,\,t\in\R,\quad |\int_s^t\tau(\xi^{\pm}(\mu))d\mu|\leq C\,.
\end{equation}
Moreover it holds
\begin{equation}\label{eq_taul}
\lim_{\xi_2\rightarrow0}\int_0^t\tau(\xi^{\pm}(s))ds =  \left\{\begin{array}{ll}0&\textrm{if } 0<t<t^*\\\beta&\textrm{if } t>t^*\end{array}\right.\,.
\end{equation}
In addition, since the lower band is not occupied initially and most of the band-to-band transfer occurs at the time $t^*$, it holds: 
\begin{equation}\label{eq_limdx}
\forall\, 0<t<t^*, \quad \lim_{\xi_2\rightarrow 0}\partial_{x_1}w_-(t,x_1^+(t),\xi_1^+(t))=0\,.
\end{equation}
In order to simplify the presentation, the proof of \eqref{eq_limdx} is postponed until the end of the present discussion. Consider first the case $t<t^*$ in \eqref{eq_dumo}. Using assumption $A2$ and equations \eqref{eq_taub} and \eqref{eq_limdx}, the dominated convergence theorem can be applied to the second term at the r.h.s. of equation \eqref{eq_dumo}, which leads to:
\begin{equation}\label{eq_lst}
\int_0^te^{-\int_s^t\tau(\xi^+(\mu))d\mu M}F(s)ds \underset{\xi_2\rightarrow0}{\longrightarrow}0\,.
\end{equation}
Combining assumption $A3$ with equations \eqref{eq_taul} and \eqref{eq_lst}, 
one can take the limit $\xi_2\rightarrow0$ in \eqref{eq_dumo} and obtain:
\begin{equation}\label{eq_trd0}
\left(\begin{array}{c}w^0_+(t,x_I+t,\xi_I-\alpha t)\\w^0_-(t,x_I+t,\xi_I-\alpha t)\end{array}\right)=\left(\begin{array}{c}w_{I,+}(x_I,\xi_I)\\w_{I,-}(x_I,\xi_I)\end{array}\right)\,.
\end{equation}
To obtain \eqref{eq_trd0}, we have used the limit below:
\begin{equation*}
\lim_{\xi_2\rightarrow0}(x_1^+(t),\xi_1^+(t))=\left\{\begin{array}{ll}(x_I+t,\xi_I-\alpha t)&\textrm{if } 0<t<t^*\\(x_I+2t^*-t,\xi_I-\alpha t)&\textrm{if } t>t^*\end{array}\right.
\end{equation*}
which follows from the explicit formula of the classical trajectories solution to \eqref{eq_caract2D}. Then, by passing to the limit $t\rightarrow t^*$ in \eqref{eq_trd0}, it follows from assumption $A4$:
\begin{equation}\label{eq_trd1}
\left(\begin{array}{c}w^r_+(t^*,x_1^*)\\w^r_-(t^*,x_1^*)\end{array}\right)=\left(\begin{array}{c}w_{I,+}(x_I,\xi_I)\\w_{I,-}(x_I,\xi_I)\end{array}\right)
\end{equation}
where $x_1^*=x_I+t^*$. The case $t>t^*$ follows the same line with the difference that \eqref{eq_lst} has to be replaced by:
\begin{equation}\label{eq_intF2I}
\int_0^te^{-\int_s^t\tau(\xi^+(\mu))d\mu M}F(s)ds = \int_0^{t^*}e^{-\int_s^t\tau(\xi^+(\mu))d\mu M}F(s)ds + \int_{t^*}^{t}e^{-\int_s^t\tau(\xi^+(\mu))d\mu M}F(s)ds\,. 
\end{equation}
The arguments used to deduce \eqref{eq_lst} can be applied here to show that the first integral at the r.h.s of equation \eqref{eq_intF2I} tends to $0$ when $\xi_2\rightarrow0$. Moreover, using assumption $A2$ and equation \eqref{eq_taub}, the second integral at the r.h.s of equation \eqref{eq_intF2I} is bounded by the quantity $C(t-t^*)$, which tends to $0$ when $t\rightarrow t^*$ (here $C$ is a constant which does not depend on $\xi_2$). Then, by taking successively the limit $\xi_2\rightarrow0$ and $t\rightarrow t^*$ in \eqref{eq_dumo}, we obtain:
\begin{equation}\label{eq_trg}
\left(\begin{array}{c}w^l_+(t^*,x_1^*)\\w^l_-(t^*,x_1^*)\end{array}\right)=e^{-\beta M}\left(\begin{array}{c}w_{I,+}(x_I,\xi_I)\\w_{I,-}(x_I,\xi_I)\end{array}\right)\,.
\end{equation}
Putting together \eqref{eq_trd1} and \eqref{eq_trg}, we get
\begin{equation}\label{eq_trans}
\left(\begin{array}{c}w^l_{+}(t^*,x_1^*)\\w^l_{-}(t^*,x_1^*)\end{array}\right)=\left(\begin{array}{cc} 1-T & T\\T & 1-T\end{array}\right)\left(\begin{array}{c}w^r_{+}(t^*,x_1^*)\\w^r_{-}(t^*,x_1^*)\end{array}\right) 
\end{equation} 
where the transition probability $T$ is given by:
\begin{equation}\label{eq_tprobm}
T = \frac{1-e^{-2\beta}}{2}
\end{equation}
and the constant $\beta$ is defined in \eqref{eq_c}. The system \eqref{eq_trans} has the form of the solution of the well known Landau-Zener problem (see \cite{Ze}). Since $x_I$ and $\xi_I>0$ are arbitrary, equation \eqref{eq_trans} is true for any $t^*>0$ and $x_1^*\in\R$.

We can now give the proof of \eqref{eq_limdx}. Differentiating \eqref{eq_siMo} with respect to $x_1$ gives the following PDE for $\partial_{x_1}w_-$:
\begin{equation}\label{eq_edpdxwm}
\partial_t\left(\partial_{x_1}w_-\right) - \frac{\xi_1}{\vert\xi\vert}\partial_{x_1}\left(\partial_{x_1}w_-\right)-\alpha\partial_{\xi_1}\left(\partial_{x_1}w_-\right)=\tau(\xi)(\partial_{x_1}w_{+}-\left(\partial_{x_1}w_-\right))\,.
\end{equation}
Consider $\varphi_s^{-}(\tilde{x}_I,\xi_I)=(x_1^-(s),\xi_1^-(s))$,  the lower band classical trajectory solution to \eqref{eq_caract2D} with an initial condition $(\tilde{x}_I,\xi_I)$ such that $\varphi_t^{-}(\tilde{x}_I,\xi_I)=\varphi_t^{+}(x_I,\xi_I)$. Then it holds
\[
\partial_{x_1}w_-(t,x_1^-(t),\xi_1^-(t))=\partial_{x_1}w_-(t,x_1^+(t),\xi_1^+(t))
\]
and it is enough to show that $\partial_{x_1}w_-(t,x_1^-(t),\xi_1^-(t))$ tends to $0$ when $\xi_2$ tends to $0$. Now, if one differentiates the map $s\mapsto \partial_{x_1}w_{-}(s,x_1^-(s),\xi_1^-(s))$, one gets using \eqref{eq_edpdxwm}: 
\[
\frac{d}{ds}\left(\partial_{x_1}w_{-}(s,x_1^-(s),\xi_1^-(s))\right)=\tau(\xi^-(s))\left(\partial_{x_1}w_{+}(s,x_1^-(s),\xi_1^-(s))-\partial_{x_1}w_-(s,x_1^-(s),\xi_1^-(s))\right)\,.
\]
By integrating the previous equation with the Duhamel formula, it follows:
\begin{equation}\label{eq_dudxwm}
\partial_{x_1}w_{-}(t,x_1^-(t),\xi_1^-(t))=\int_0^te^{-\int_s^t\tau(\xi^-(\mu))d\mu}\tau(\xi^-(s))\partial_{x_1}w_{+}(s,x_1^-(s),\xi_1^-(s))ds
\end{equation}
where we have used that $w_{I,-}=0$. Using assumption $A2$, together with the equations \eqref{eq_taub} and \eqref{eq_taul}, we can conclude that $\forall\, 0<t<t^*$, the integral at the r.h.s. of \eqref{eq_dudxwm} tends to $0$ when $\xi_2\rightarrow 0$.
\begin{remark}
Although the system \eqref{eq_trans} has the correct form, the formula \eqref{eq_tprobm} for the transition probability is different from the Landau-Zener transition probability \eqref{eq_tprob} obtained in section \ref{sec_heur} and which writes for our particular potential:
\begin{equation}\label{eq_transr}
T = e^{-\frac{\pi \xi_2^2}{h|\alpha|}}\,.
\end{equation}
It will be verified numerically in section \ref{sec_Monum} that the band transmission corresponding to the effective model \eqref{eq_siMo} is given by \eqref{eq_tprobm} whereas the band transmission corresponding to the model \eqref{eq_liouvMo} is given by the correct transition probability \eqref{eq_transr}. This signifies that the approximation made in \cite{MoSu2} to obtain \eqref{eq_siMo} from \eqref{eq_liouvMo} is not correct.
\end{remark}
\begin{remark}
If the coefficient $\tau$ is replaced by $\frac{1}{\xi_2}$, \eqref{eq_siMo} becomes a hyperbolic relaxation system (see \cite{Ji}, \cite{Pa}) and the solution to \eqref{eq_siMo} satisfies $w_+=w_-$ when $\xi_2\rightarrow 0$. Since $\int w_+ +\int w_-=1$ this leads to:
\begin{equation}\label{eq_wdem}
\int w_+ =\int w_-=\frac{1}{2}\,.
\end{equation}
We will observe numerically in Figure \ref{fig_popEf} that \eqref{eq_wdem} is true when the time is big enough.   
\end{remark}

\section{Numerical results}\label{sec_numres}

In this section, the results provided by the surface hopping algorithm are compared with the reference level populations given by \eqref{eq_lpop} where the solution $u^h(t)$ is computed numerically using an accurate method to solve the Dirac equation \eqref{eq_dir}. In particular, it is verified numerically in section \ref{sec_quantnum} that the spectral method is more accurate than the finite difference method. The image of the operator $\Pi_{\pm}(hD)$ is computed by using discrete Fourier transform (DFT) and Fourier multiplication. A comparison of the surface hopping algorithm with the models \eqref{eq_liouvMo} and \eqref{eq_siMo} is given in section \ref{sec_Monum}.

\subsection{Setup and quantum level simulations}\label{sec_quantnum}

We suppose that the initial data $u^h_I$ is such that its Fourier transform $\hat{u}^h_I$ satisfies the relation:
\begin{equation}\label{eq_Fci}
\hat{u}^h_I(\xi) = \hat{f}^h(\xi)\chi_+(h\xi)
\end{equation}
where $\chi_+(\xi)$ is defined in \eqref{eq_vectprop} and $\hat{f}^h(\xi)$ is the Fourier transform of some function $f^h\in L^2(\R^2)$. Such an initial data satisfies the non-interference condition \eqref{eq_condID}. We remark that, like for the level population \eqref{eq_lpop}, $\chi_+(h\xi)$ can be replaced by $\chi_+(\xi)$ in \eqref{eq_Fci}.

If $f^h$ is bounded in $L^2(\R^2)$, using Lemma 2.3 in \cite{GMMP} again, we obtain the following approximation for the scalar Wigner transform of $u^h_I$:
\begin{equation}\label{eq_swui}
w^h\left[u^h_I\right] = W^h\left[f^h\right] + o(1)
\end{equation}
in $\mathcal{D}'\left(\R^2_x\times(\R^2_\xi\setminus\{\xi=0\})\right)$ when $h\rightarrow 0$, where $W^h\left[f^h\right]=w^h(f^h,f^h)$. By plugging \eqref{eq_swui} in \eqref{eq_dtwt}, we obtain the following asymptotics for the initial value of the diagonal terms of the Wigner matrix:
\begin{equation}\label{eq_diagwigap}
w_{I,+}^h = W^h\left[f^h\right] + o(1)\,, \quad w_{I,-}^h = o(1)
\end{equation}
in $\mathcal{D}'\left(\R^2_x\times(\R^2_\xi\setminus\{\xi=0\})\right)$ when $h\rightarrow 0$. In the present section, the initial upper level function is an $h$-scaled Gaussian wave packet:
\[
f^h(x) = (\pi h)^{-\frac{1}{2}}e^{-\frac{|x-x_0^h|^2}{2h}+i\frac{\xi_0.(x-x_0^h)}{h}}
\]
with center $x_0^h\in\R^2$, momentum $\xi_0\in\R^2$ and norm $\Vert f^h\Vert_{L^2(\R^2)}=1$. Its $h$-scaled Fourier transform and its Wigner transform can be computed explicitly. Indeed, 
\[
\mathcal{F}^hf^h(\xi) = (\pi h)^{-\frac{1}{2}}e^{-\frac{|\xi-\xi_0|^2}{2h}-i\frac{x_0^h.\xi}{h}}
\]
where $\forall u\in L^2(\R^2),\,\mathcal{F}^hu(\xi)=h^{-1}\mathcal{F}u(h^{-1}\xi)$. Moreover, 
\[
W^h\left[f^h\right](x,\xi) = (\pi h)^{-2}e^{-\frac{|x-x_0^h|^2}{h}-\frac{|\xi-\xi_0|^2}{h}}\,.
\]
In this section, the parameter $h$ is equal to its physical value given by \eqref{eq_defh}. The effective mass of the electron is given by
\[
m=0.067m_e
\]
where $m_e$ is the mass of the electron, the Fermi velocity is taken as $v_F = 10^6\,m.s^{-1}$ and, having in mind the simulation of devices of size equal to hundreds of nanometers as in \cite{MoSu2,HPA}, we will take the reference length equal to $L=500\,nm$. Then, the numerical value of the parameter $h$ corresponding to formula \eqref{eq_defh} is:
\[
h = 3.4557\times10^{-3}
\]
which is small enough for the problem to be considered in the semiclassical regime. The simulation domain $\Omega$ is equal to
\[
\Omega = [-10\sqrt{h},10\sqrt{h}]\times[-5\sqrt{h},5\sqrt{h}]\,.
\]
In the next section, we solve the Dirac equation \eqref{eq_dir} for different choice of the potential $U$ by using the time-splitting spectral method (TSSM) presented in Appendix \ref{sec_Dsol}.

\subsubsection{The Klein tunneling}\label{sec_kt}
 
\begin{figure}
\begin{center}
\includegraphics[width=\linewidth]{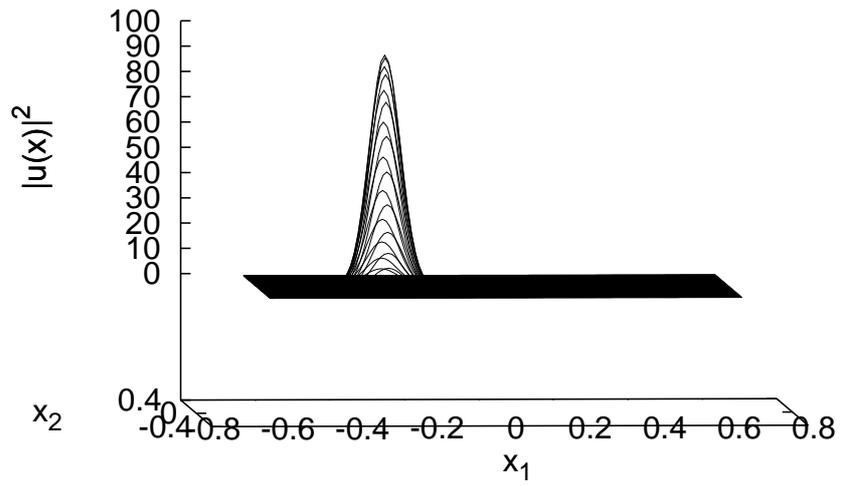}\\
\includegraphics[width=\linewidth]{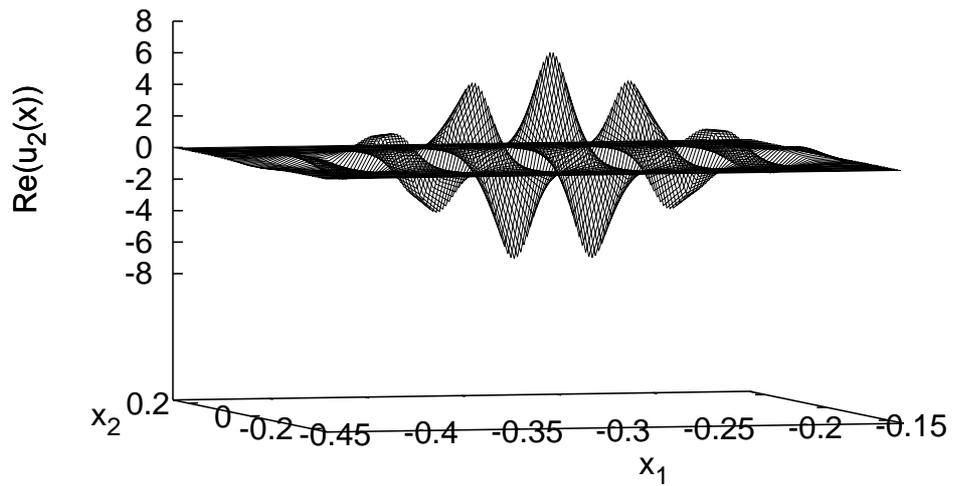}
\caption{Initial condition for the TSSM: $|u_I(x)|^2$ (top) and, in a smaller region, $\textrm{Re}\left((u_I)_2(x)\right)$ (bottom).}\label{fig_condI}
\end{center}
\end{figure}

The potential is equal to $U=v_01_{x_1\geq0}$ which corresponds to a Klein step. For this potential, we compare the TSSM  and the Finite difference time domain method (FDTD) in \cite{HPA}. The center and the momentum of the initial Gaussian wave packet $f^h$ are chosen as
\[
x_0^h = (-5\sqrt{h},0), \quad \xi_0 =\frac{1}{2}(1,0)\,.
\]
The height of the Klein step and the simulation time are respectively
\[
v_0 = 2|\xi_0| \, \textrm{ and } \, t_f = 13\sqrt{h}\,.
\]
For the TSSM, the number of grid points are given by
\[
N_1 = 126, \quad N_2 = 66 
\]
and for the FDTD by
\[
N_1 = 850, \quad N_2 = 426 \,.
\]
For the two methods, the time step size is
\[
\Delta t = \frac{\Delta x_1}{\sqrt{2}}\,,
\]
where $\Delta x_1$ and $\Delta x_2$ are the mesh sizes in the $x_1$ and $x_2$ 
directions respectively. In the present section and in the following section, the supscript $h$ will be omitted for the solution $u^h(t,x)$ to \eqref{eq_dir} and for the initial condition $u^h_I(x)$. Using the discrete Fourier 
transform (DFT) \eqref{eq_defDFT}, resp. its inverse, to approximate the Fourier transform, resp. the inverse Fourier transform, we get the following approximation of the projectors $\Pi_{\pm}(hD)u(t^n,x_j)$: 
\begin{equation}\label{eq_projn}
\Pi_{\pm}(hD)u^n_j = \frac{1}{N_1N_2}\sum_{k\in\mathcal{K}}\Pi_{\pm}(h\xi_k)\widehat{\left(u^n\right)}_ke^{i\xi_k.(x_j-a)}, \quad j\in\mathcal{J}
\end{equation}
where $\widehat{\left(u^n\right)}_k$ is the Fourier coefficient defined in \eqref{eq_defDFT} and the discretization is the same as in Appendix \ref{sec_Dsol}. Then, using formula \eqref{eq_lpop}, the  approximation $P^h_{dir,\pm}(t^n)$ of the level populations $P^h_{\pm}(t^n)$ of the Dirac equation is given by:
\begin{equation}\label{eq_lpopdir}
P^h_{dir,\pm}(t^n) = \left\Vert\Pi_{\pm}(hD)u^n\right\Vert_2^2
\end{equation}
where $\Pi_{\pm}(hD)u^n_j$ is defined by \eqref{eq_projn} and for $u=(u_j)_{j\in\mathcal{J}}$, $u_j\in\mathbb{C}^2$, we have
\begin{equation}\label{eq_defnorm}
\left\Vert u\right\Vert_2^2 = \sum_{j\in\mathcal{J}}|u_j|^2\Delta x_1\Delta x_2\,.
\end{equation}
Similarly, the initial condition $u_I(x_j)$ defined by \eqref{eq_Fci} is approximated by the formula:
\[
\left(u_I\right)_j = \frac{2\pi}{(b_1-a_1)(b_2-a_2)}\sum_{k\in\mathcal{K}}\hat{f}^h(\xi_k)\chi_+(h\xi_k)e^{i\xi_k.x_j}, \quad j\in\mathcal{J}\,.
\]   
We remark that in the previous formula, we used the exact value of the Fourier transform $\hat{f}^h(\xi_k)$ instead of the DFT $\widehat{\left(f^h\right)}_k$. The initial data for the TSSM with the parameters described above is represented in Figure \ref{fig_condI}.

\begin{figure}
\begin{center}
\includegraphics[width=\linewidth]{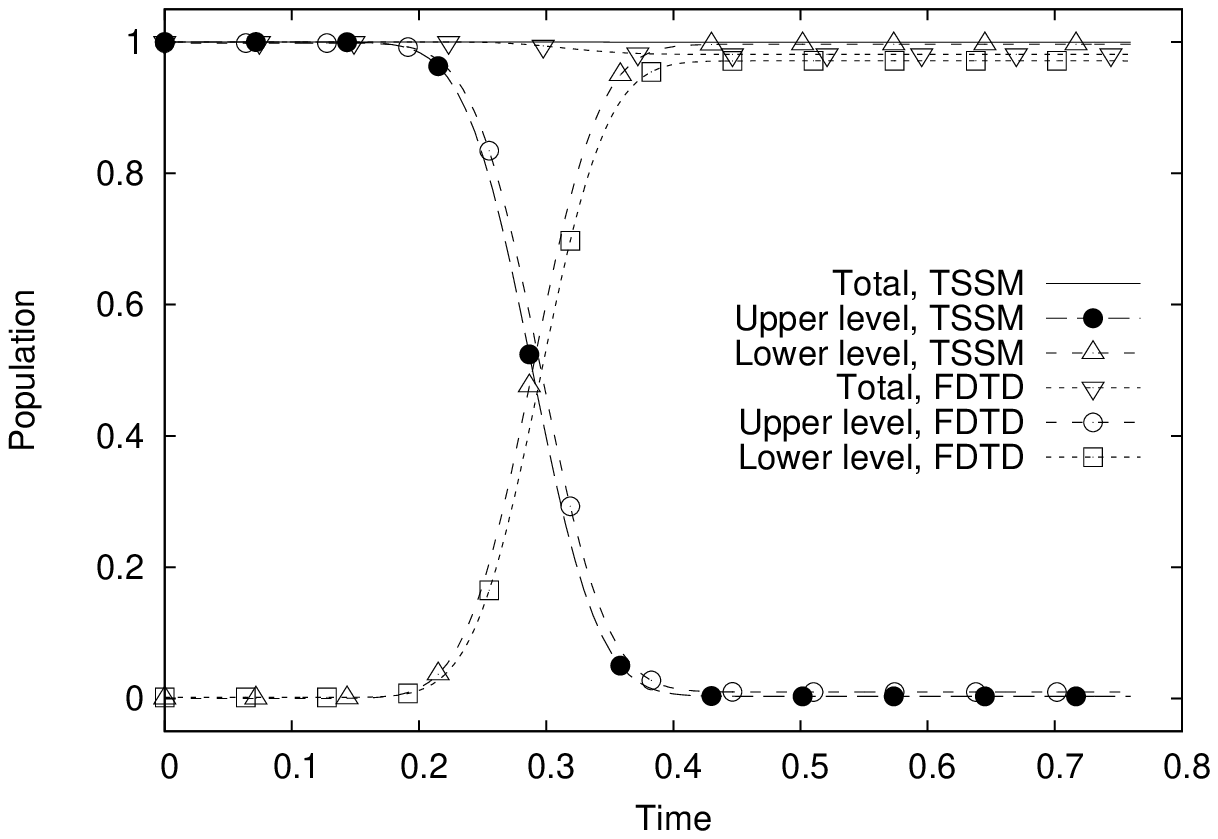}
\caption{Time evolution of the level populations with a Klein step potential: comparison of the TSSM and the FDTD.}\label{fig_levpop}
\end{center}
\end{figure}

 In Figure \ref{fig_levpop}, we depict the evolution with respect to the time $t^n$ of the level populations $P^h_{dir,\pm}(t^n)$ provided by the TSSM and by the FDTD. The curve with title Upper level, resp. Lower level, refers to the plus sign, resp. minus sign, in \eqref{eq_lpopdir} and the curve with title Total refers to the total population $\left\Vert u^n\right\Vert_2^2$. Initially, the charge is carried completely by the upper level, then non-adiabatic transfer occurs at time $t=0.28$ and the charge is almost all transferred to the lower level. Moreover, we remark that the TSSM is more accurate than the FDTD. Indeed, the first method conserves  the total mass  \eqref{eq_mconsd} whereas for the second method the total mass decreases at the hopping time (it was shown in \cite{HPA} that the quantity conserved by the FDTD is not the mass but a related functional). To reduce this mass loss, the number of spatial points has to be chosen big enough which increases the CPU time of the method. In particular, the CPU times corresponding to the simulations of Figure \ref{fig_levpop} are the following: $493.5268\,s$ for the FDTD, $1.3961\,s$ for the TSSM.

\begin{figure}
\begin{center}
\includegraphics[width=\linewidth]{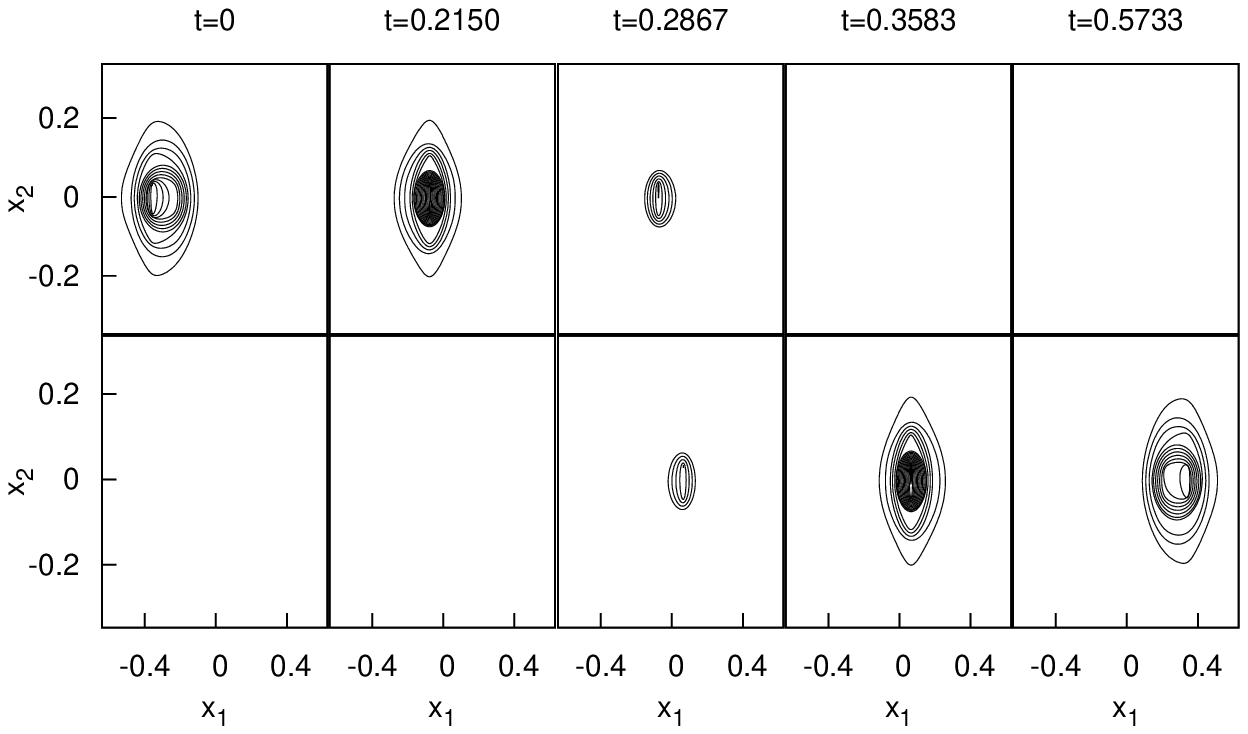}
\caption{For different times and with respect to the space variable: representation of the position density of the projection of the wave function $u$ solution to \eqref{eq_dir} where $U=v_01_{x_1\geq0}$ and the initial condition is given by \eqref{eq_Fci}. The upper half corresponds to $|\Pi_+(hD)u^n_j|^2$ and the lower half to $|\Pi_-(hD)u^n_j|^2$ where the projectors are computed using \eqref{eq_projn}. The solution $u$ is computed with the TSSM.}\label{fig_projdens}
\end{center}
\end{figure}

 We remark in Figure \ref{fig_projdens} that the transition occurs when the charge reaches the Klein step which is a classically forbidden region for the Upper level. The band transition to the Lower level makes the Klein step an allowed region for the particles, which can tunnel in the region $x_1\geq 0$ with a probability almost equal to $1$. This is the Klein paradox, see \cite{Th}. 

\subsubsection{Case $U=\alpha x_1$, $\alpha>0$}

\begin{figure}
\begin{center}
\includegraphics[width=0.9\linewidth]{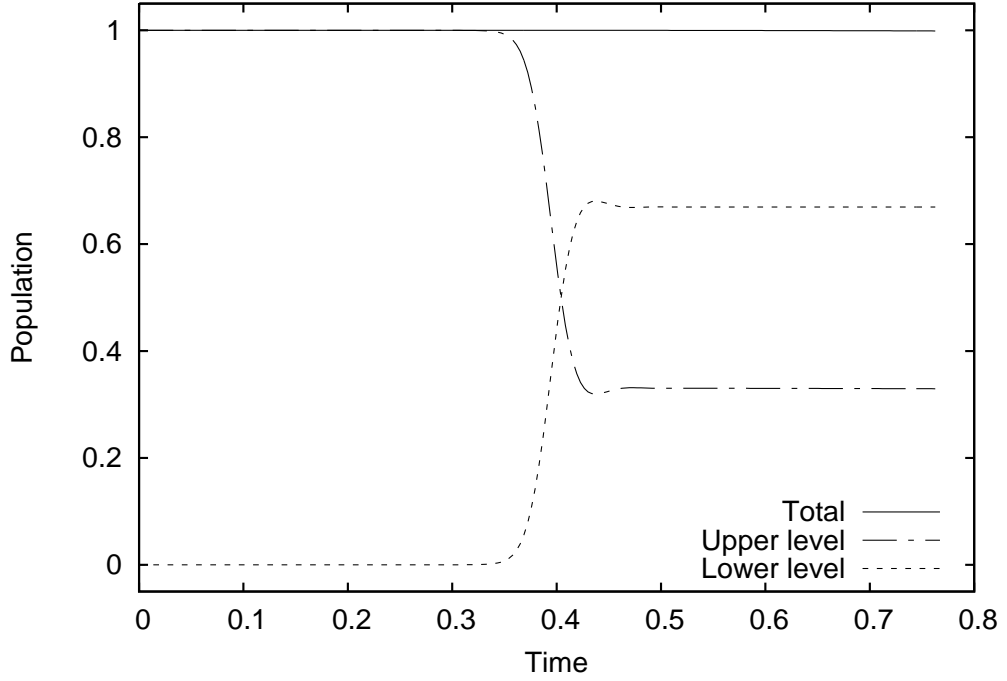}
\caption{Time evolution of the level populations for the potential $U=\alpha x_1$, $\alpha>0$.}\label{fig_popPL}
\end{center}
\end{figure}

In this section, the potential is given by $U=\alpha x_1$, $\alpha>0$. The Dirac equation \eqref{eq_dir} is solved using the TSSM. The center and the momentum of the initial Gaussian wave packet $f^h$ are taken to be
\[
x_0^h = (-5\sqrt{h},0), \quad \xi_0 = (1,0)\,.
\]
The level position density of the solution $u(t,x)$ to \eqref{eq_dir} is supported around $(x^{\pm}(t),\xi^{\pm}(t))$, solution to:
\[
\displaystyle \left\{\begin{array}{ll}\frac{d}{dt}x^{\pm}(t)= \pm\frac{\xi^{\pm}(t)}{\vert\xi^{\pm}(t)\vert}, & x^{\pm}(0)=x_0^h\\
\frac{d}{dt}\xi^{\pm}(t)=-\alpha\left(1,0\right), & \xi^{\pm}(0)=\xi_0
\end{array}\right.\,.
\] 
The solution of the above problem can be computed explicitly. It is given by:
\begin{equation}\label{eq_flowexp}
x^{\pm}(t) = x^* \mp\left(|t-t^*|,0\right), \quad \xi^{\pm}(t)=\xi_0-\alpha\left(t,0\right)
\end{equation}
where $t^*=\frac{(\xi_0)_1}{\alpha}$ is the time such that $\xi^+(t^*)=0$ and the non-adiabatic transfer occurs. The point $x^*=x^+(t^*)=x_0^h+(t^*,0)$ is the point where the hopping occurs. The coefficient $\alpha$ is chosen such that the potential at the hopping point is equal to $U(x^*)=v_0$ where $v_0=\frac{|\xi_0|}{4}$. This leads to:
\[
\alpha = \frac{v_0-(\xi_0)_1}{(x_0^h)_1}\,.
\]
The simulation stops at time:
\[
t_f = 13\sqrt{h}
\]
and the number of space points are given by
\[
N_1 = 160, \quad N_2 = 80\,. 
\]
The time step size is equal to
\[
\Delta t = \frac{\Delta x_1}{\sqrt{2}}\,.
\]
The initial data $\left(u_I\right)_j$, the Lower level and Upper level populations $P^h_{dir,\pm}(t^n)$ and the Total population are computed as explained in section \ref{sec_kt}. The time evolution of the level populations is represented in Figure \ref{fig_popPL}. The numerical hopping is observed around the time $t=0.39$ which is accurate enough compared to the predicted value $t^*=0.3919$ given by \eqref{eq_flowexp}. Contrary to the case of the Klein step potential, a significant part of the charge stays on the upper level.

\begin{figure}
\begin{center}
\includegraphics[width=\linewidth]{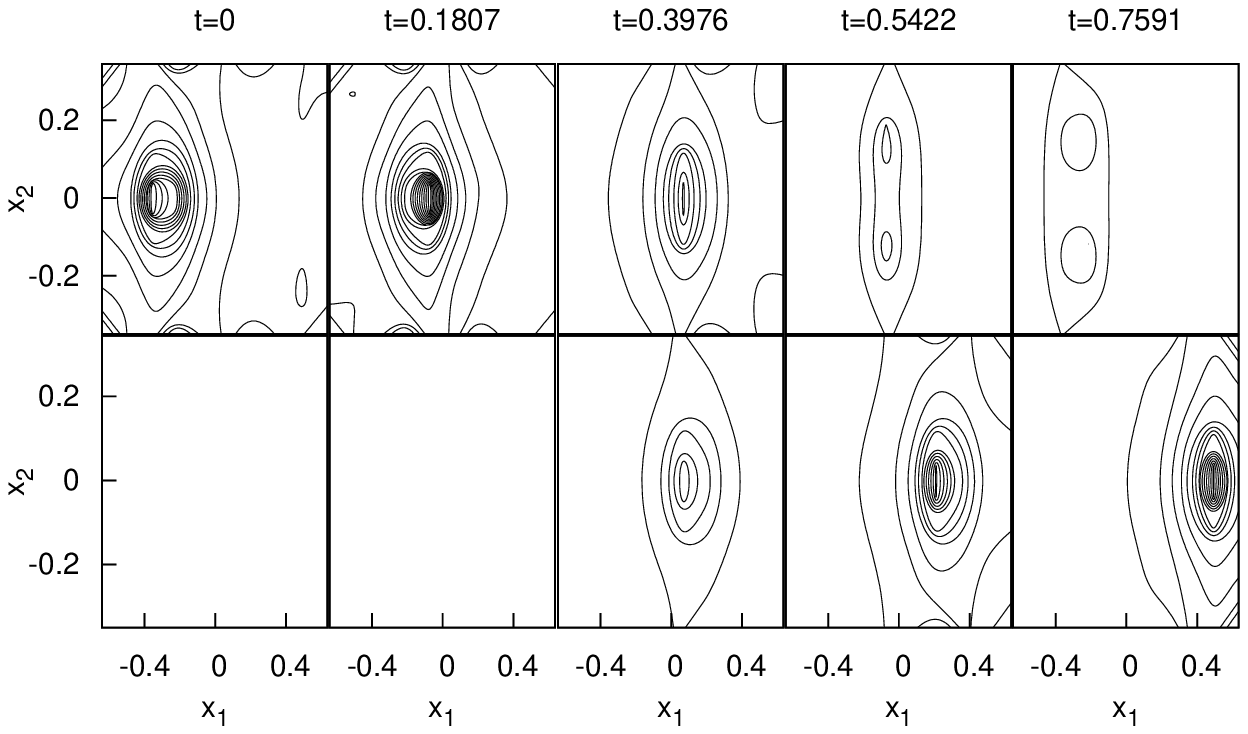}
\caption{For different times and with respect to the space variable: representation of the position density of the projection of the wave function $u$ solution to \eqref{eq_dir} where $U=\alpha x_1$ and the initial condition is given by \eqref{eq_Fci}. The upper half corresponds to $|\Pi_+(hD)u^n_j|^2$ and the lower half to $|\Pi_-(hD)u^n_j|^2$ where the projectors are computed using \eqref{eq_projn}. The solution $u$ is computed with the TSSM.}\label{fig_projdenspl}
\end{center}
\end{figure}

The predicted value of the position corresponding to the hopping is $x^*=(9.7976\times10^{-2},0)$. We remark in Figure \ref{fig_projdenspl} that the transition occurs when the charge reaches $x^*$. As for the Klein step potential, the band transition to the lower level allows the particles to tunnel in the region $x_1\geq x^*_1$. However, for the potential considered here, we can observe that the part remaining on the upper level is reflected. This could have been predicted from equation \eqref{eq_flowexp}. Indeed, for $t\geq t^*$, the classical flow corresponding to the upper level (plus sign) moves to the left with respect to $x^*$ whereas the classical flow corresponding to the lower level (minus sign) moves to the right with respect to $x^*$.  

\subsection{The surface hopping algorithm}\label{sec_ressh}

In the present section, the potential $U$ is equal to $U=\alpha x_1$ where $\alpha = 15$. The initial condition is as described in section \ref{sec_quantnum}. The center and the momentum of the initial Gaussian wave packet $f^h$ are taken as
\begin{equation}\label{eq_centsh}
x_0^h = (-5\sqrt{h},0), \quad \xi_0 =(1,0)\,.
\end{equation}
The diagonal terms $w_{\pm}^h(t,x,\xi)$ of the Wigner matrix defined by \eqref{eq_wh} are computed using the surface hopping algorithm presented in section \ref{sec_alg}. Then, for different values of the parameter $h$ and of the simulation time $t_f$, the surface hopping level populations $P^h_{sh,\pm}(t_f)$ are computed from $w_{\pm}^h(t_f,.,.)$ by using \eqref{eq_lpopsh}. The results provided by the surface hopping algorithm are compared to the level populations $P^h_{dir,\pm}(t_f)$ of the Dirac equation computed as explained in section \ref{sec_kt}.

For the level populations $P^h_{dir,\pm}(t_f)$, the Dirac equation \eqref{eq_dir} is solved using the TSSM in the simulation domain
\[
\Omega = [-11\sqrt{h},11\sqrt{h}]\times[-5\sqrt{h},5\sqrt{h}]\,.
\] 
Using equation \eqref{eq_mconsd}, we can choose a time step size which depends only on the simulation time $t_f$ as follows: 
\begin{equation}\label{eq_dtref}
\Delta t = t_f/250\,.
\end{equation}
The number of space points are given by: for $h=10^{-1}$ and $h=10^{-2}$  
\[
N_1 = 250, \quad N_2 = 126\,, 
\]
for $h=10^{-3}$ 
\[
N_1 = 300, \quad N_2 = 150\,,
\]
for $h=10^{-4}$ 
\[
N_1 = 850, \quad N_2 = 426\,. 
\]
For the level populations $P^h_{sh,\pm}(t_f)$, the initial sampling is chosen as follows. We consider a uniform $J\times J$ discretization of the domain of the $x$ variable:
\[
[a_1,b_1]\times[a_2,b_2]
\]
where $a_1=a_2=x_0^h-5\sqrt{h}$ and $b_1=b_2=x_0^h+5\sqrt{h}$ and a uniform $K\times K$ discretization of the $\xi$ variable domain
\[
[c_1,d_1]\times[c_2,d_2]
\]
where $c_1=c_2=\xi_0-5\sqrt{h}$ and $d_1=d_2=\xi_0+5\sqrt{h}$. The grid points $x_j, \, j=1,\hdots,J^2$ and $\xi_k, \, k=1,\hdots,K^2$ are ordered such that
\[
|f^h(x_1)|^2\geq \hdots \geq |f^h(x_{J^2})|^2\,,
\]
\[
|\mathcal{F}^hf^h(\xi_1)|^2\geq \hdots \geq |\mathcal{F}^hf^h(\xi_{K^2})|^2\,.
\]
Then, we determine the minimal numbers $N_x$ of points of the $x$ variable and $N_{\xi}$ of the $\xi$ variable such that
\[
\sum_{j=1}^{N_x}|f^h(x_j)|^2 \Delta x^2 \geq 1 - \textrm{tol}_x, \quad
\sum_{k=1}^{N_{\xi}}|\mathcal{F}^hf^h(\xi_k)|^2 \Delta \xi^2 \geq 1 - \textrm{tol}_{\xi}
\]
where $\Delta x=\frac{b_1-a_1}{J}=\frac{b_2-a_2}{J}$, $\Delta \xi=\frac{d_1-c_1}{K}=\frac{d_2-c_2}{K}$ and $\textrm{tol}_x$, $\textrm{tol}_{\xi}$ are well chosen tolerances. The phase space points $(x_j,\xi_k)$, $j=1,\hdots,N_x$, $k=1,\hdots,N_{\xi}$ are ordered such that:
\[
W^h\left[f^h\right](x_1,\xi_1)\geq \hdots \geq W^h\left[f^h\right](x_{N_xN_{\xi}},\xi_{N_xN_{\xi}})
\]   
and we determine the minimal integer $N$ such that:
\[
\sum_{k=1}^{N}W^h\left[f^h\right](x_k,\xi_k)\Delta x^2\Delta \xi^2 \geq 1 - \textrm{tol}
\]
where $\textrm{tol}$ is a well chosen tolerance. This gives the first step of the algorithm of section \ref{sec_alg}. Indeed, the initial sampling is the set:
\[
\{(x_k,\xi_k,+)\in\R^2_x\times\R^2_\xi\times\{-,+\}\,; \quad k=1,...,N\}
\]
where, using equation \eqref{eq_diagwigap}, the associated weights $w_k\in\R$ can be approximated by:
\[
w_k = W^h\left[f^h\right](x_k,\xi_k)\,.
\]
In the present case, the hopping surface $S$ defined by \eqref{eq_minsurf} is equal to:
\[
S = \{ (x,\xi)\in\R^4;\, \xi_1=0  \}
\]
and the second order derivative of the function $s\mapsto\vert\xi^{\pm}(s)\vert^2$, defined by the characteristics solution to \eqref{eq_caract}, is equal to:
\[
\left(\vert\xi^{\pm}(s)\vert^2\right)''=2\alpha^2>0\,.
\]
Therefore, for the potential considered here, the points of extremal gap are all minimas.\\
Moreover, the classical flow $(x^{\pm}(t),\xi^{\pm}(t))$ solution to \eqref{eq_caract} is given by
\begin{equation}\label{eq_trajUa}
x^{\pm}(t) =  x \pm \int_0^t\frac{\xi^{\pm}(s)}{|\xi^{\pm}(s)|}ds, \quad \xi^{\pm}(t) =  \xi-\alpha t(1,0)\,.
\end{equation} 
In equation \eqref{eq_trajUa}, the formula of the momentum is explicit and the hopping transport step can be simplified. Indeed, for $k=1,\hdots,N$, if $0<\frac{(\xi_k)_1}{\alpha}<t_f$, the trajectory $(x_k(t),\xi_k(t))$, $0<t<t_f$ defined by \eqref{eq_decev} will pass through a point $(x_k^*,\xi_k^*)\in S$ at the time $t^*=\frac{(\xi_k)_1}{\alpha}$ and non-adiabatic transfer occurs. For such a $k$, the weight is changed such that for $t>t^*$
\[
w_k(t) = \left(1-T^*\right)w_k(t^*) 
\]
and a new particle is created on the lower band with index $l>N$. For the new particle, the associated weight is such that for $t>t^*$
\[
w_l(t) = T^*w_k(t^*)\,. 
\]
In the above equations, the transition rate $T^*$ is equal to
\[
T^*=T(x_k^*,\xi_k^*)
\] 
where $T(x,\xi)$ is given by \eqref{eq_transr}. The surface hopping level populations defined by \eqref{eq_lpopsh} are then given by:
\[
P^h_{sh,+}(t_f) = \sum_{k=1}^{N}w_k(t_f)\Delta x^2\Delta \xi^2\,, \quad P^h_{sh,-}(t_f) = \sum_{l=N+1}^{M}w_l(t_f)\Delta x^2\Delta \xi^2\,.
\]
For all the tests, the size of the sampling grids are equal to:
\[
J=K=16\,.
\]
The sampling tolerances are taken as:
\[
\textrm{tol} = 10^{-6}, \quad \textrm{tol}_x = \textrm{tol}_{\xi} = 10^{-3}\times\textrm{tol}\,.
\]
For such parameters, the number of particles obtained numerically for the initial sampling is equal to $N=6981$.

\begin{figure}
\begin{center}
\includegraphics[width=0.9\linewidth]{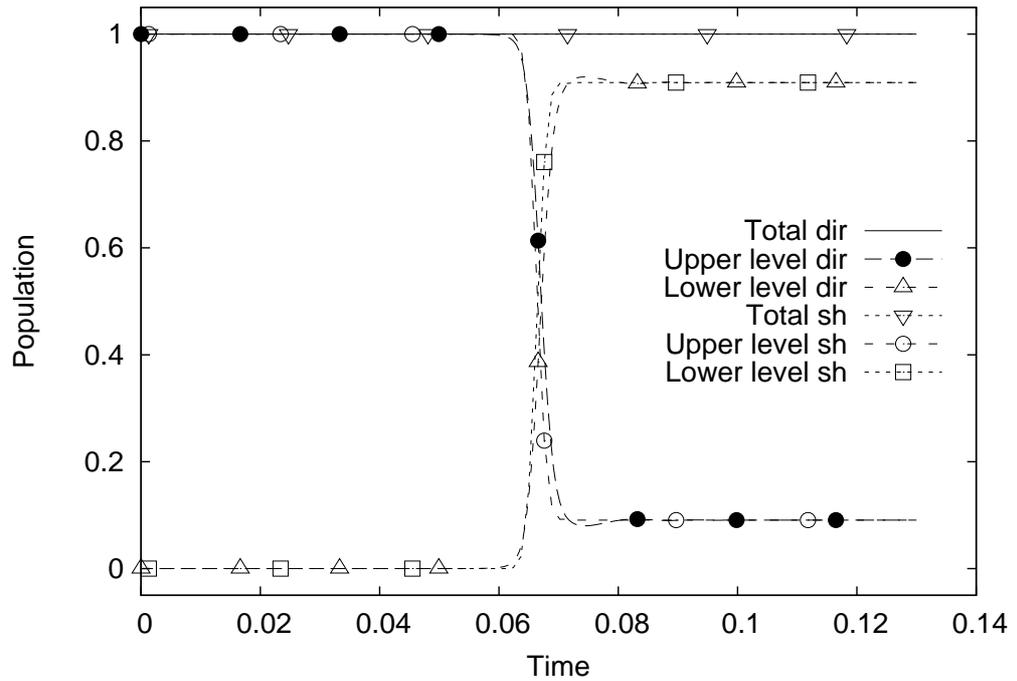}
\caption{For $h=10^{-3}$ and $U=\alpha x_1$, $\alpha=15$: time evolution of the level populations provided by the Dirac solver and the surface hopping algorithm.}\label{fig_popDS}
\end{center}
\end{figure}

For the simulation time $t_f = 0.13$ and for different values of the semiclassical parameter $h$, the level populations obtained by the two methods are listed in Table \ref{tab_levp}. We remark that the reference lower level population $P^h_{dir,-}(t_f)$ increases when $h$ increases. This is due to the fact that for larger values of $h$, the transition process is slower and, at the time $t_f$, the post-transition relaxation observed in Figure \ref{fig_popDS}, has no yet happened. We notice that, for $h$ smaller or equal to $10^{-2}$, the surface hopping level populations $P^h_{sh,\pm}(t_f)$ are almost constant with respect to $h$. This can be explained by the fact that, for the potential considered, the transition rate depends only on the variable $\xi_2/\sqrt{h}$ which does not depend on $h$ when considered on the sampling points. The column CPU dir, resp. CPU sh, denotes the CPU time required by the Dirac solver, resp. the surface hopping method, to complete the simulations. The time required to choose the initial sampling is not taken into account in CPU sh. It appears that the surface hopping method is very interesting. Indeed, when $h$ decreases, the numerical cost of the Dirac solver increases whereas the surface hopping method provides very accurate results for an almost constant and much smaller CPU time.
\begin{table}
\begin{center}
\begin{tabular}{|c|c|c|c|c|c|c|}
\hline 
 h & $P^h_{dir,+}(t_f)$ & $P^h_{sh,+}(t_f)$ & $P^h_{dir,-}(t_f)$ & $P^h_{sh,-}(t_f)$ & CPU dir ($s$) & CPU sh ($s$) \\
\hline
 $10^{-1}$ & $7.48058\times 10^{-2}$ & $9.07411\times 10^{-2}$ & $0.9251938$ & $0.9092579$ & $10.1206$ & $0.1400$ \\
\hline
 $10^{-2}$ & $8.98146\times 10^{-2}$ & $9.06980\times 10^{-2}$ & $0.9101854$ &  $0.9093010$ & $10.1766$  & $0.1360$ \\
\hline
 $10^{-3}$ & $9.06998\times 10^{-2}$ & $9.06980\times 10^{-2}$ & $0.9093002$ & $0.9093010$ & $14.7849$ & $0.1320$ \\
\hline
 $10^{-4}$ & $9.06978\times 10^{-2}$ & $9.06980\times 10^{-2}$ & $0.9093015$ & $0.9093010$ & $198.3884$ & $0.1320$ \\
\hline
\end{tabular}
\caption{At time $t_f = 0.13$ and for different values of the semiclassical parameter $h$: level populations obtained by the Dirac solver and the surface hopping method.}\label{tab_levp}
\end{center}
\end{table}

\begin{figure}
\begin{center}
\includegraphics[width=0.9\linewidth]{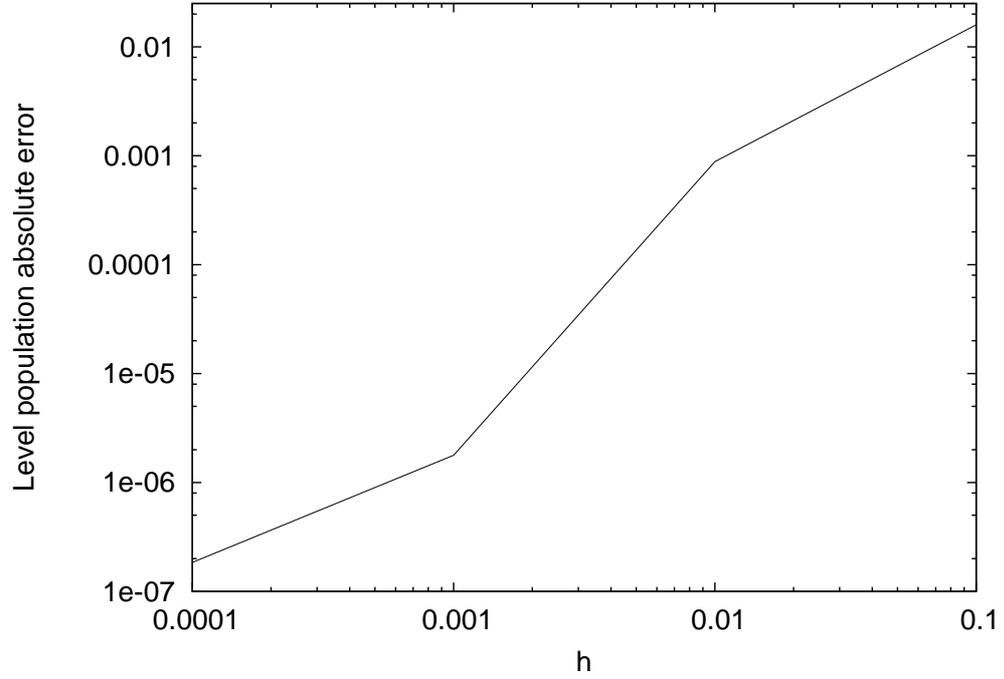}
\caption{At time $t_f = 0.13$, logarithmic plot of the absolute error \eqref{eq_levpoper} of the level population corresponding to the upper level while varying the semiclassical parameter $h=10^{-p}$, $p=1,2,3,4$.}\label{fig_levpoper}
\end{center}
\end{figure}

We depict in Figure \ref{fig_levpoper} the absolute error of the level population corresponding to the upper level:
\begin{equation}\label{eq_levpoper}
\mathcal{E}^h_+(t_f) := |P^h_{dir,+}(t_f)-P^h_{sh,+}(t_f)|
\end{equation}
with respect to the semiclassical parameter $h$. We remark that the surface hopping level populations converge numerically when $h\rightarrow0$ to the level populations obtained at the quantum level.

For $h=10^{-3}$, the time evolution of the level populations provided by the Dirac solver and the surface hopping algorithm are depicted in Figure \ref{fig_popDS}. The curve with title Upper level dir, resp. Lower level dir, refers to $P^h_{dir,+}(t)$, resp. $P^h_{dir,-}(t)$. The same is true when dir is replaced by sh. The total population of the Dirac equation, curve titled Total dir, is obtained as explained in section \ref{sec_kt} and the surface hopping total population, curve titled Total sh, corresponds to $P^h_{sh,+}(t)+P^h_{sh,-}(t)$. We observe that the numerically obtained surface hopping total population is conserved. For a smaller CPU time, the time evolution of the level populations provided by the surface hopping algorithm agrees well with the one provided by the Dirac solver. Indeed, for the two methods, the charge is initially carried completely by the upper level, then non-adiabatic transfer occurs at time $t=\frac{(\xi_0)_1}{\alpha}=\frac{1}{15}=0.0667$ (time required by the classical flow to reach the crossing set $\{\xi=0\}$ starting from the momentum $\xi_0$ of the initial Gaussian wave packet) and the great majority of the charge is transferred to the lower level. The CPU time required to complete the simulation is $14.7849\,s$ for the Dirac solver and $5.6803\,s$ for the surface hopping method. For the second method, the CPU time is bigger than the time indicated in Table \ref{tab_levp}. This is due to the fact that the surface hopping curves in Figure \ref{fig_popDS} are obtained by repeating the surface hopping algorithm described in the present section for the sequence of times $t_f=0.13n/100$, $1\leq n\leq100$.  

\subsection{Simulations using the models of \cite{MoSu2}}\label{sec_Monum}

In this section we will compare numerically the models \eqref{eq_liouvMo} and \eqref{eq_siMo} with the two dimensional version of our surface hopping algorithm. The comparison will be performed with $h=10^{-3}$, $t_f=0.13$ and with the same potential as is section \ref{sec_ressh}. The initial condition is the 2D version of the initial condition in section \ref{sec_ressh}. In other words:
\[
w^h_{+}|_{t=0}= W^h[f^h], \quad w^h_{-}|_{t=0}=0 
\]
is replaced by
\begin{equation}\label{eq_ciMo}
w_{+}(0,x_1,\xi_1)= (\pi h)^{-1}e^{-\frac{(x_1-(x_0^h)_1)^2}{h}-\frac{(\xi_1-(\xi_0)_1)^2}{h}}, \quad w_{-}(0,x_1,\xi_1)=0\,.
\end{equation}
The functions $w_{\pm}$ and $w_{i}$ depend on $\xi_2$. However, since $\xi_2$ plays the role of a parameter, this dependence is not written on the l.h.s. of the equations in \eqref{eq_ciMo}. The center and the momentum of the initial Gaussian wave packet are given by \eqref{eq_centsh}. 

\subsubsection{The two dimensional surface hopping algorithm}\label{sec_sh2d}

The $x_2$-independent solutions to the surface hopping algorithm presented in section \ref{sec_alg} are obtained by replacing the system \eqref{eq_liouv} with its two dimensional version:
\begin{equation}\label{eq_liouv2d}
\partial_t w_{\pm} \pm \frac{\xi_1}{\vert\xi\vert}\partial_{x_1}w_{\pm}-\alpha\partial_{\xi_1}w_{\pm}=0, \quad \xi_1\neq 0\,.
\end{equation}
More precisely, equation \eqref{eq_liouv2d} is used for the time evolution of $w_{\pm}$ as long as the classical trajectories $\varphi^{\pm}_t(x_1,\xi_1)$ solution to \eqref{eq_caract2D} are away from the hopping surface $\{\xi_1=0\}$. Then, the hopping transport is described as follows.\\
Consider an initial set of sampling points
\[
\{(x_{1,k},\xi_{1,k},+)\in\R_{x_1}\times\R_{\xi_1}\times\{-,+\}\,; \quad k=1,...,N\}
\]
with associated weights $w_k\in\R$ given by:
\[
w_k = w_+(0,x_{1,k},\xi_{1,k})\,.
\]
For $t\geq 0$ small enough:
\[
(x_{1,k}(t),\xi_{1,k}(t)) = \varphi^+_t(x_{1,k},\xi_{1,k}), \quad w_k(t) = w_k 
\]
where, using \eqref{eq_caract2D}, we have $\xi_{1,k}(t)=\xi_{1,k}-\alpha t$. Then, for $k=1,\hdots,N$, if $0<\frac{\xi_{1,k}}{\alpha}<t_f$, the classical trajectory is such that $\xi_{1,k}(t^*)=0$ at the time $t^*=\frac{\xi_{1,k}}{\alpha}$ and non-adiabatic transfer occurs. For such a $k$, the weight is changed such that for $t>t^*$
\[
w_k(t) = \left(1-T^*\right)w_k(t^*) 
\]
and a new particle is created on the lower band with index $l>N$. For the new particle, the associated weight is such that for $t>t^*$
\[
w_l(t) = T^*w_k(t^*)\,. 
\]
In the above equations, the transition rate $T^*$ is equal to:
\[
T^*=e^{-\frac{\pi \xi_2^2}{h|\alpha|}}\,.
\] 
The surface hopping level populations at the final time $t_f$ are then given by:
\begin{equation}\label{eq_popsh2d}
P^h_{sh,+}(t_f) = \sum_{k=1}^{N}w_k(t_f)\Delta x_1\Delta \xi_1\,, \quad P^h_{sh,-}(t_f) = \sum_{l=N+1}^{M}w_l(t_f)\Delta x_1\Delta \xi_1
\end{equation}
where $\Delta x_1$ and $\Delta \xi_1$ are the mesh sizes in the $x_1$ and
$\xi_1$ directions respectively.

\subsubsection{Comparison with the models in \cite{MoSu2}}

For equations \eqref{eq_liouvMo} and \eqref{eq_siMo}, the simulation domain is:
\[
(x_1,\xi_1) \in \left[-11\sqrt{h},11\sqrt{h}\right]\times\left[(\xi_0)_1-\alpha t_f-5\sqrt{h},(\xi_0)_1+5\sqrt{h}\right]\,.
\]
They are solved with periodic boundary conditions on a uniform grid with $500$ points in the $x_1$-direction and $500$ points in the $\xi_1$-direction. The time step size is chosen such that the condition of stability of the upwind method is satisfied \cite{Le}. More precisely, we take:
\[
\Delta t = \left(\frac{1}{\Delta x_1}+\frac{\alpha}{\Delta\xi_1}\right)^{-1}\,.
\] 
Equations \eqref{eq_liouvMo} and \eqref{eq_siMo} are solved using a time-splitting method where the free equation (without source term) is solved using a dimensional splitting: the free problem is splitted in two one-dimensional problems and each one dimensional problem is solved using the one-dimensional second order upwind MUSCL scheme (see\cite{Le}). The source term is integrated in time using a RK2 method for equation \eqref{eq_liouvMo} (in order to preserve the second order accuracy) and exactly for equation \eqref{eq_siMo}.

For equation \eqref{eq_liouvMo}, we take:
\[
w_{i}|_{t=0}=0
\]
in addition to the initial condition \eqref{eq_ciMo}. We remark that, since the MUSCL method is written for real valued functions, the third equation in \eqref{eq_liouvMo} has to be splitted in two equations, one for the real part of $w_{i}$ and one for its imaginary part. 

The level populations provided by the asymptotic model are defined by:
\[
P^h_{am,\pm}(t) = \int_{\R^2} w_{\pm}(t,x_1,\xi_1)dx_1d\xi_1
\]
where $w_{\pm}$ is the solution to \eqref{eq_liouvMo} and the level populations provided by the effective model are defined by:
\[
P^h_{em,\pm}(t) = \int_{\R^2} w_{\pm}(t,x_1,\xi_1)dx_1d\xi_1
\]
where $w_{\pm}$ is the solution to \eqref{eq_siMo}. For the level populations $P^h_{sh,\pm}(t_f)$ defined in \eqref{eq_popsh2d}, the initial sampling is chosen as follows. We consider a uniform $J\times K$ discretization of the domain of the $(x_1,\xi_1)$ variables:
\[
[a,b]\times[c,d]
\]
where $a=(x_0^h)_1-5\sqrt{h}$, $b=(x_0^h)_1+5\sqrt{h}$ and $c=(\xi_0)_1-5\sqrt{h}$, $d=(\xi_0)_1+5\sqrt{h}$. The grid points $(x_{1,j},\xi_{1,k})$, $j=1,\hdots,J$, $k=1,\hdots,K$ are ordered such that:
\[
w_+(0,x_{1,1},\xi_{1,1}) \geq \hdots \geq w_+(0,x_{1,J\times K},\xi_{1,J\times K})
\]   
and we determine the minimal integer $N$ such that:
\[
\sum_{k=1}^{N}w_+(0,x_{1,k},\xi_{1,k})\Delta x_1\Delta \xi_1 \geq 1 - \textrm{tol}
\]
where $\Delta x_1=\frac{b-a}{J}$, $\Delta \xi_1=\frac{d-c}{K}$ and $\textrm{tol}$ is a well chosen tolerance. Then, the initial sampling is the set:
\[
\{(x_{1,k},\xi_{1,k},+)\in\R_{x_1}\times\R_{\xi_1}\times\{-,+\}\,; \quad k=1,...,N\}\,.
\]
For all the tests in this section, the number of grid points for the two dimensional surface hopping algorithm is:
\[
J=100, \quad K=100
\]
and the tolerance is
\[
\textrm{tol}=10^{-9}\,.
\]

\begin{figure}
\begin{center}
\includegraphics[width=0.9\linewidth]{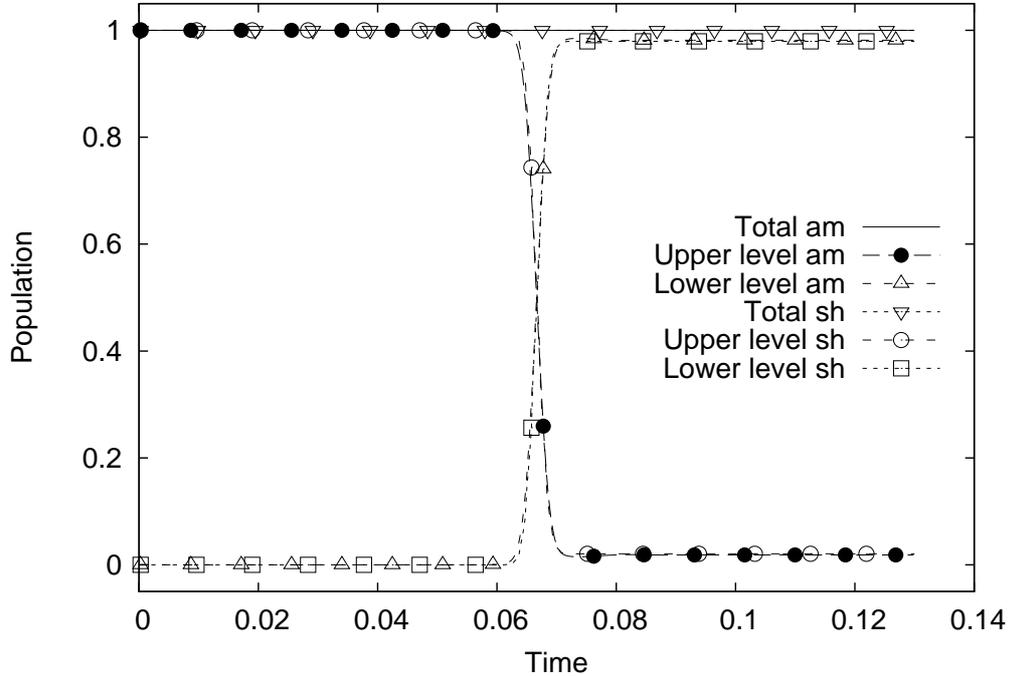}
\caption{For $h=10^{-3}$, $\xi_2=10^{-2}$ and $U=\alpha x_1$, $\alpha=15$: time evolution of the level populations provided by the asymptotic model \eqref{eq_liouvMo} and the two dimensional surface hopping algorithm presented in section \ref{sec_sh2d}.}\label{fig_popMo}
\end{center}
\end{figure}

We depict in Figure \ref{fig_popMo} the time evolution of the level populations provided by the asymptotic model \eqref{eq_liouvMo} and the two dimensional surface hopping algorithm. The curve with title Upper level am, resp. Lower level am, refers to $P^h_{am,+}(t)$, resp. $P^h_{am,-}(t)$. The total population, curve titled Total am, corresponds to $P^h_{am,+}(t)+P^h_{am,-}(t)$. The same is true when am is replaced by sh. The surface hopping curves are obtained by repeating the surface hopping algorithm described in section \ref{sec_sh2d} for the sequence of times $t_f=0.13n/500$, $1\leq n\leq500$. The behavior is the same as in section \ref{sec_ressh}: the charge is initially carried completely by the upper level, then hopping occurs at time $t=0.0667$ and the great majority of the charge is transferred to the lower level. The time evolution of the level populations provided by the asymptotic model fits well the one provided by the two dimensional surface hopping: the surface hopping algorithm validates the model \eqref{eq_liouvMo}.

\begin{figure}
\begin{center}
\includegraphics[width=0.9\linewidth]{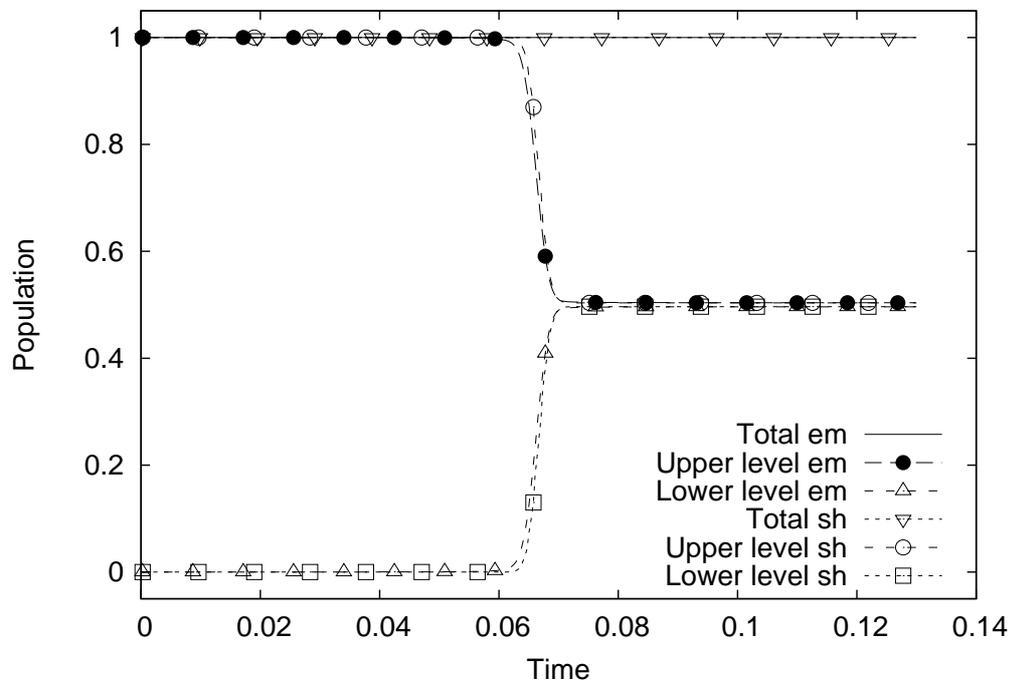}
\caption{For $h=10^{-3}$, $\xi_2=10^{-2}$ and $U=\alpha x_1$, $\alpha=15$: time evolution of the level populations provided by the effective model \eqref{eq_siMo} and the two dimensional surface hopping algorithm presented in section \ref{sec_sh2d} where the transition probability $T^*$ is replaced by \eqref{eq_tprobm}\eqref{eq_c}.}\label{fig_popEf}
\end{center}
\end{figure}

Figure \ref{fig_popEf} shows that non-adiabatic transfer is only partially recovered by  model \eqref{eq_siMo}. We depict the time evolution of the level populations provided by the asymptotic model \eqref{eq_siMo} and the two dimensional surface hopping algorithm where the transition probability $T^*$ is replaced by \eqref{eq_tprobm}\eqref{eq_c}. The curve with title Upper level em, resp. Lower level em, refers to $P^h_{em,+}(t)$, resp. $P^h_{em,-}(t)$. The total population, curve titled Total em, corresponds to $P^h_{em,+}(t)+P^h_{em,-}(t)$. The same is true when em is replaced by sh. As in Figure \ref{fig_popMo}, hopping occurs at time $t=0.0667$. However, using the simplified model \eqref{eq_siMo},  only about a half of the charge is transferred to the lower level which is close to the transition probability given by \eqref{eq_tprobm}\eqref{eq_c} and substantially different from Figure \ref{fig_popMo}.

\begin{figure}
\begin{center}
\includegraphics[width=0.9\linewidth]{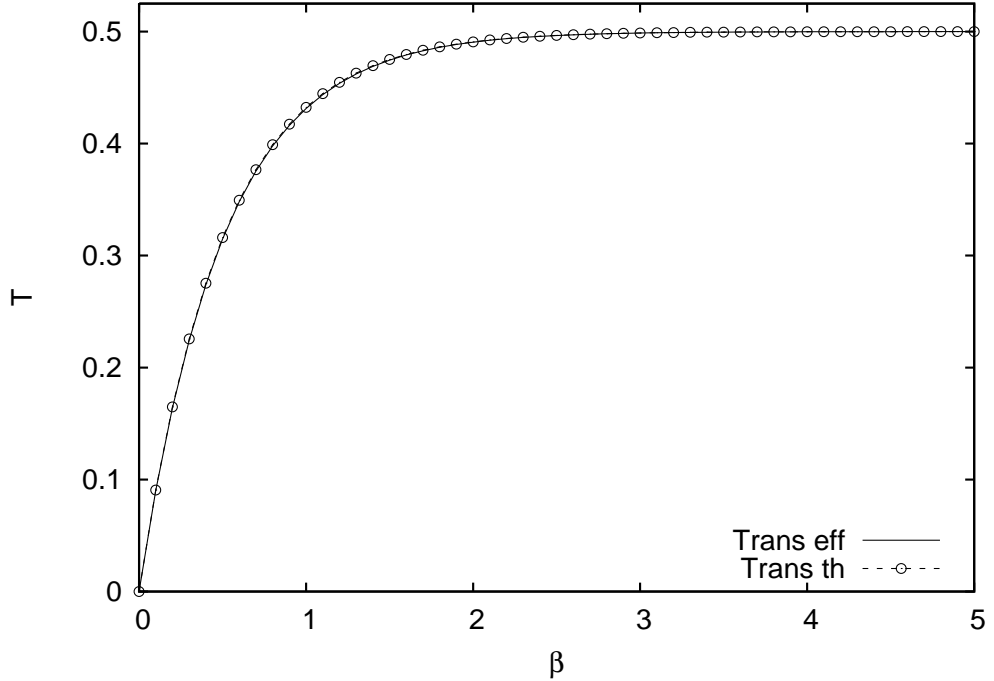}
\caption{For $h=10^{-3}$, $\xi_2=10^{-2}$, $t_f=0.13$ and $U=\alpha x_1$, $\alpha=15$: numerical verification of the formula \eqref{eq_tprobm} for the transition probability corresponding to the model \eqref{eq_siMo}.}\label{fig_verT}
\end{center}
\end{figure}

In Figure \ref{fig_verT}, we verify numerically that the transition probability corresponding to the effective model \eqref{eq_siMo} is given by the formula \eqref{eq_tprobm}. The curve with title Trans eff refers to $P^h_{em,-}(t_f)$, the population provided by the effective model \eqref{eq_siMo} on the lower level and at the final time, for $31$ different values of the constant $\beta$ distributed uniformly on the interval $[0,5]$. The constant $\beta$ given by \eqref{eq_c} is made arbitrary in \eqref{eq_tprobm} by replacing the coefficient $\tau$ appearing in \eqref{eq_siMo} by $\tilde{\tau}=\frac{\beta}{\pi^2/4}\tau$. The curve with title Trans th is the representation of the coefficient $T$ given by \eqref{eq_tprobm} for the same values of $\beta$. We remark that the two curves are very close which validates the limit $\xi_2\rightarrow0$ performed in section \ref{sec_thmo} for the solution to \eqref{eq_siMo}.

\begin{appendices}

\section{Time-splitting spectral method (TSSM)}\label{sec_Dsol}

The solution $u(t,x)$  to \eqref{eq_dirad} is computed on the domain:
\[
\Omega = [a_1,b_1]\times[a_2,b_2]
\]
using a time-splitting spectral method as in \cite{HJMSZ}. For $r=1,2$, we choose the spatial mesh size $\Delta x_r=\frac{b_r-a_r}{N_r}$ in the $r$ direction for a given integer $N_r$. Define a uniform grid
\begin{equation}\label{eq_defma}
x_j = a + (j_1\Delta x_1,j_2\Delta x_2), \quad j\in\mathcal{J}
\end{equation}
where 
\[
\mathcal{J} = \{ j=(j_1,j_2)\in\mathbb{N}^2 \, | \, 0 \leq j_1 < N_1\,,\, 0 \leq j_2 < N_2 \}
\]
and $a=(a_1,a_2)$. For a given time step size $\Delta t$, let $u^n_j$ denote the numerical approximation of $u(t^n,x_j)$ at the time $t^n=n\Delta t$, $n\geq0$. Then, $u^{n+1}_j$ is computed from  $u^n_j$ by decomposing the problem \eqref{eq_dirad} in the two sub-problems
\begin{equation}\label{eq_dirli_r}
ih\partial_tu = \left[-ih\sigma_1\partial_{x_1}-ih\sigma_2\partial_{x_2}\right]u
\end{equation}
and
\begin{equation}\label{eq_dynpot}
ih\partial_tu = Uu\,.
\end{equation}
The free Dirac equation \eqref{eq_dirli_r} is solved using a spectral method in space and exact time integration, whereas equation \eqref{eq_dynpot} can be integrated exactly on $[t^n,t^{n+1}]$. To discretize \eqref{eq_dirli_r}, we introduce the trigonometric interpolant of $u$:
\begin{equation}\label{eq_interp}
\tilde{u}(t,x) = \frac{1}{N_1N_2}\sum_{k\in\mathcal{K}}\widehat{u(t)}_ke^{i\xi_k.(x-a)} 
\end{equation}
where $\tilde{u}(t,x_j)=u(t,x_j)$. In equation \eqref{eq_interp}, we have
\[
\mathcal{K} = \{ k=(k_1,k_2)\in\mathbb{Z}^2 \, | \, -\frac{N_1}{2} \leq k_1 < \frac{N_1}{2}\,,\,-\frac{N_2}{2} \leq k_2 < \frac{N_2}{2} \}
\]  
and
\[
\xi_k = 2\pi\left(\frac{k_1}{b_1-a_1},\frac{k_2}{b_2-a_2}\right)\,.
\]
For a given function $f:\R^2\rightarrow \C^2$, $\hat{f}_k$ is the discrete Fourier transform (DFT) of $f$ defined by
\begin{equation}\label{eq_defDFT}
\hat{f}_k = \sum_{j\in\mathcal{J}}f_je^{-i\xi_k.(x_j-a)}
\end{equation}
where $f_j=f(x_j)$. When $f=(f_j)_{j\in\mathcal{J}}$ is a sequence with $f_j\in\C^2$, the DFT of $f$ is the sequence $\hat{f}_k$ defined by equation \eqref{eq_defDFT}. Inserting \eqref{eq_interp} into \eqref{eq_dirli_r}, and using the orthogonality of the set $\left\{e^{i\xi_k.(x-a)},\,k\in\mathcal{K}\right\}$ with respect to the scalar product of $L^2(\Omega)$, one gets the following system of ODE
\begin{equation*}
\frac{d}{dt}\widehat{u(t)}_k = -iB(\xi_k)\widehat{u(t)}_k\,.
\end{equation*}
For any $t_0\in\R$, the exact solution to the previous equation is given by:
\[
\widehat{u(t)}_k = M(t-t_0,\xi_k)\widehat{u(t_0)}_k
\] 
where
\[
M(\delta,\xi) = e^{-i\delta B(\xi)}\,.
\]
Applying an inverse discrete Fourier transform, one obtains the following expression for the approximation $u(t)_j$ of $u(t,x_j)$:
\[
u(t)_j = \frac{1}{N_1N_2}\sum_{k\in\mathcal{K}}M(t-t_0,\xi_k)\widehat{u(t_0)}_ke^{i\xi_k.(x_j-a)}\,.
\]
We remark that using the eigenvectors \eqref{eq_vectprop}, the matrix $B(\xi)$ can be diagonalized as follows:
\[
B(\xi) = P(\xi)D(\xi)P(\xi)^*
\]
where
\[
D(\xi) = \textrm{diag}\left(|\xi|,-|\xi|\right)\,, \ \quad P(\xi) = \left(\chi_+(\xi),\chi_-(\xi)\right) 
\]
and therefore
\begin{equation}\label{eq_matM}
M(\delta,\xi) = P(\xi)e^{-i\delta D(\xi)}P(\xi)^* = \left(\begin{array}{cc} \cos(\delta|\xi|)&-i\sin(\delta|\xi|)(\xi_1-i\xi_2)/|\xi|\\-i\sin(\delta|\xi|)(\xi_1+i\xi_2)/|\xi|&\cos(\delta|\xi|)\end{array}\right)\,.
\end{equation}
The Strang splitting is the second order method constructed as follows: solve the first subproblem \eqref{eq_dirli_r} over only a half time step of length $\frac{\Delta t}{2}$. Then, we use the result as data for a full time step on the second subproblem \eqref{eq_dynpot} and finally take another half time step on \eqref{eq_dirli_r}. This leads to the following method:
\begin{align}
& u^*_j = \frac{1}{N_1N_2}\sum_{k\in\mathcal{K}}M(\Delta t/2,\xi_k)\widehat{\left(u^n\right)}_ke^{i\xi_k.(x_j-a)}\nonumber\\
& u^{**}_j = e^{-\frac{i}{h}U_j\Delta t}u^*_j\label{eq_TSSM}\\
& u^{n+1}_j = \frac{1}{N_1N_2}\sum_{k\in\mathcal{K}}M(\Delta t/2,\xi_k)\widehat{\left(u^{**}\right)}_ke^{i\xi_k.(x_j-a)}\,.\nonumber
\end{align}
\begin{remark} 
It follows directly from the unitarity of the DFT and of the matrix $M(\delta,\xi)$ given by \eqref{eq_matM} that the TSSM \eqref{eq_TSSM} conserves the discrete total charge, i.e.:
\begin{equation}\label{eq_mconsd}
\Vert u^n \Vert_2 = \Vert u^0 \Vert_2\,, \quad \forall n\geq0 
\end{equation}
where the norm $\Vert.\Vert_2$ is defined in \eqref{eq_defnorm}. Moreover, when the potential is equal to a constant $U=U_0$, the Dirac equation \eqref{eq_dirad} admits the following plane wave solutions:
\[
u(t,x) = \chi_{\pm}(k)e^{\frac{i}{h}\left(k.x-(U_0\pm|k|)t\right)}\,,
\]
which are integrated exactly by the TSSM \eqref{eq_TSSM} if $k\in\R^2$ 
satisfies $\frac{k}{h}=\xi_{k'}$ for some $k'\in\mathcal{K}$. In the case of the Schr\"odinger equation, an analysis of the stability of the TSSM was performed in \cite{FGL} for initial conditions which are close to plane waves.
\end{remark}
\begin{remark}
In the applications, the solution of the Dirac equation \eqref{eq_dirad} moves away from the simulation box. Therefore, in addition to the periodic boundary conditions provided by the spectral method, we use absorbing boundary layers at the edge of the simulation box, see e.g. \cite{RKPK}.
\end{remark}
\end{appendices}

\begin{description}
\item[Acknowledgements] A.F. acknowledges Clotilde Fermanian Kammerer for 
 discussion about section \ref{sec_rig}. 
\end{description}

\end{document}